\documentclass[11pt]{article}
\usepackage{amssymb,amsfonts}
\begin{document}

\begin{center}
\Large Apparent singularities of Fuchsian equations, \\
and the Painlev\'e VI equation and Garnier systems
\end{center}
\medskip

\centerline{ R.\,R.\,Gontsov\footnote{Institute for Information
Transmission Problems, Moscow, Russia, rgontsov@inbox.ru.} and
I.V. Vyugin\footnote{Institute for Information Transmission
Problems, Moscow, Russia, ilya\_vyugin@mail.ru} }
\medskip

\begin{abstract}
We study movable singularities of Garnier systems (and Painlev\'e
VI equations) using the connection of the latter with
isomonodromic deformations of Fuchsian systems. Questions on the
existence of solutions for some inverse monodromy problems are
also considered.
\end{abstract}
\bigskip

{\bf \S1. Introduction}
\medskip

In the middle of the XIXth century B.\,Riemann \cite{Ri}
considered the problem of the construction of a linear
differential equation
\begin{eqnarray}\label{urav}
\frac{d^pu}{dz^p}+b_1(z)\frac{d^{p-1}u}{dz^{p-1}}+\ldots+b_p(z)u=0
\end{eqnarray}
with the prescribed regular singularities
$a_1,\ldots,a_n\in\overline{\mathbb C}$ (which are the poles of
the coefficients) and prescribed monodromy.

Recall that a singular point $a_i$ of the equation (\ref{urav}) is
said to be {\it regular} if any solution of the equation is of no
more than a polynomial (with respect to $1/|z-a_i|$) growth near
$a_i$.

By L.\,Fuchs's theorem (\cite{Fu}, see also \cite{Ha}, Th. 12.1) a
singular point $a_i$ is regular if and only if the coefficient
$b_j(z)$ has at this point a pole of order $j$ or lower
$(j=1,\ldots,p)$. Linear differential equations with regular
singular points only are called {\it Fuchsian}.

The monodromy of a linear differential equation describes a
branching pattern of its solutions near singular points. It is
defined as follows. In a neighbourhood of a non-singular point
$z_0$ we consider a basis $(u_1,\ldots,u_p)$ in the solution space
of the equation (\ref{urav}). Analytic continuations of the
functions $u_1(z),\ldots,u_p(z)$ along an arbitrary loop $\gamma$
outgoing from $z_0$ and lying in $\overline{\mathbb C}\setminus
\{a_1,\ldots,a_n\}$ transform the basis $(u_1,\ldots,u_p)$ into a
(in general case another) basis $(\tilde u_1,\ldots,\tilde u_p)$.
The two bases are related by means of a non-singular transition
matrix $G_{\gamma}$ corresponding to the loop $\gamma$:
$$
(u_1,\ldots,u_p)=(\tilde u_1,\ldots,\tilde u_p)G_{\gamma}.
$$
The map $[\gamma]\mapsto G_{\gamma}$ (which depends only on the
homotopy class $[\gamma]$ of the loop $\gamma$) defines the
representation
$$
\chi: \pi_1(\overline{\mathbb C}\setminus\{a_1,\ldots,a_n\},z_0)
\longrightarrow{\rm GL}(p,\mathbb C)
$$
of the fundamental group of the space $\overline{\mathbb
C}\setminus \{a_1,\ldots, a_n\}$ in the space of non-singular
complex matrices of size $p$. This representation is called the
{\it monodromy} of the equation (\ref{urav}).

By the {\it monodromy matrix} of the equation (\ref{urav}) at a
singular point $a_i$ (with respect to the basis
$(u_1,\ldots,u_p)$) we mean the matrix $G_i$ corresponding to a
simple loop $\gamma_i$ encircling $a_i$, so that
$G_i=\chi([\gamma_i])$. The matrices $G_1,\ldots,G_n$ are the
generators of the {\it monodromy group} of the equation
(\ref{urav}). They satisfy the relation $G_1\ldots G_n=I$ implied
by the condition $\gamma_1\ldots\gamma_n=e$ for the generators of
the fundamental group (here and below $I$ is the identity matrix).

If one considers another basis
$(u'_1,\ldots,u'_p)=(u_1,\ldots,u_p)C$, $C\in{\rm GL}(p,\mathbb
C)$, of the solution space of the equation (\ref{urav}), then the
corresponding monodromy matrices change as follows:
$G'_i=C^{-1}G_iC$. In a similar way the matrices $G_i$ depend on
the choice of an initial point $z_0$. So one sees that the
monodromy of a linear differential equation is defined up to a
conjugation by a constant non-singular matrix and it is more
precisely to say that the monodromy is an element of the space
$$
{\cal M}_a={\rm Hom}\left(\pi_1(\overline{\mathbb
C}\setminus\{a_1,\ldots,a_n\}), {\rm GL}(p,\mathbb C)\right) /{\rm
GL}(p,\mathbb C)
$$
of conjugacy classes of representations of the group
$\pi_1(\overline{\mathbb C}\setminus\{a_1,\ldots,a_n\})$.

A.\,Poincar\'e \cite{Po} has established that the number of
parameters determining a Fuchsian equation of order $p$ with $n$
singular points is less than the dimension of the space ${\cal
M}_a$ of monodromy representations, if $p>2, n>2$ or $p=2, n>3$
(see also \cite{AB}, pp. 158--159). Hence in the construction of a
Fuchsian equation with the given monodromy there arise (besides
$a_1,\ldots,a_n$) so-called {\it apparent} singularities at which
the coefficients of the equation have poles but the solutions are
single-valued meromorphic functions, i.~e., the monodromy matrices
at these points are identity matrices. Below by apparent singular
points of an equation we mean these very singularities. Thus, in
general case the Riemann problem has a negative solution.

A similar problem for systems of linear differential equations is
called the {\it Riemann--Hilbert problem}. This is the problem of
the construction of a {\it Fuchsian system}\footnote{Note that by
Sauvage's theorem singular points of a Fuchsian system are regular
(see \cite{Ha}, Th. 11.1).}
\begin{eqnarray}\label{syst}
\frac{dy}{dz}=\left(\sum_{i=1}^n\frac{B_i}{z-a_i}\right)y, \qquad
y(z)\in{\mathbb C}^p, \qquad B_i\in{\rm Mat}(p,{\mathbb C}),
\end{eqnarray}
with the given singularities $a_1,\ldots,a_n$ (if $\infty$ is not
a singular point of the system, then $\sum_{i=1}^nB_i=0$) and
monodromy
\begin{eqnarray}\label{repr}
\chi: \pi_1(\overline{\mathbb C}\setminus\{a_1,\ldots,a_n\},z_0)
\longrightarrow{\rm GL}(p,\mathbb C).
\end{eqnarray}

One defines the monodromy of a linear system in the same way as
for the scalar equation (\ref{urav}); one merely needs to consider
in place of a basis $(u_1,\ldots,u_p)$ in the solution space of
the equation a {\it fundamental matrix} of the system, i.~e.,
matrix whose columns form a basis in the solution space of the
system.

A counterexample to the Riemann--Hilbert problem was obtained by
A.\,A.\,Bolibrukh (\cite{Bo1}, see also \cite{AB}, Ch. 5). The
solution of this problem has a more complicated history than that
of the Riemann problem for scalar Fuchsian equations (before
A.\,A.\,Bolibrukh it had long been wrongly regarded as solved in
the affirmative; for details see \cite{AL}, \cite{Bo6}).

Alongside Fuchsian equations consider the famous non-linear
differential equations of second order
--- the {\it Painlev\'e VI equation} (${\rm P_{VI}}$) and {\it Garnier systems.}

The equation ${\rm P_{VI}}(\alpha,\beta,\gamma,\delta)$ is the
non-linear differential equation
\begin{eqnarray}\label{PVI}
\frac{d^2u}{dt^2}&=&\frac12\left(\frac1u+\frac1{u-1}+\frac1{u-t}\right)\left(\frac{du}{dt}\right)^2-
\left(\frac1t+\frac1{t-1}+\frac1{u-t}\right)\frac{du}{dt}+\nonumber \\
& &+\frac{u(u-1)(u-t)}{t^2(t-1)^2}\left(\alpha+\beta\,\frac
t{u^2}+\gamma\,\frac{t-1}{(u-1)^2}+\delta\,\frac{t(t-1)}{(u-t)^2}\right)
\end{eqnarray}
of second order with respect to the unknown function $u(t)$, where
$\alpha,\beta,\gamma,\delta$ are complex parameters. This equation
has three fixed singular points --- $0, 1, \infty$. Its movable
singularities (which depend on the initial conditions) can be
poles only. In such a case one says that an equation satisfies the
{\it Painlev\'e property}. The general ${\rm P_{VI}}$ equation
(\ref{PVI}) was first written down by R.\,Fuchs \cite{Fu1} (son of
L.\,Fuchs) and was added to the list of the equations now known as
the {\it Painlev\'e I--VI equations} by Painlev\'e's student
B.\,Gambier \cite{Ga}. Among the non-linear differential equations
of second order satisfying the Painlev\'e property, only the
equations of this list in general case can not be reduced to the
known differential equations for elementary and classical special
functions. The ${\rm P_{VI}}$ equation is the most general because
all the other ${\rm P_{I-V}}$ equations can be derived from it by
certain limit processes after the substitution of the independent
variable $t$ and parameters (see \cite{GP}, Ch. III, \S1.2).
R.\,Fuchs suggested two approaches to obtaining the ${\rm P_{VI}}$
equation. The first one (which we will discuss in this paper)
deals with isomonodromic deformations of Fuchsian systems. The
second, more geometrical, approach uses elliptic integrals.

The Garnier system ${\cal G}_n(\theta)$ depending on $n+3$ complex
parameters $\theta_1,\ldots,\theta_{n+2},\theta_{\infty}$ is a
completely integrable Hamiltonian system (see \cite{GP}, Ch. III,
\S4)
\begin{eqnarray}\label{garnier}
\frac{\partial u_i}{\partial a_j}=\frac{\partial H_j}{\partial
v_i}, \qquad \frac{\partial v_i}{\partial a_j}= -\frac{\partial
H_j}{\partial u_i}, \qquad i,j=1,\ldots,n,
\end{eqnarray}
with certain Hamiltonians $H_i=H_i(a,u,v,\theta)$ rationally
depending on $a=(a_1,\ldots,a_n)$, $u=(u_1,\ldots,u_n)$,
$v=(v_1,\ldots,v_n)$, $\theta=(\theta_1,\ldots,\theta_{n+2},
\theta_{\infty})$. In the case $n=1$ the Garnier system ${\cal
G}_1(\theta_1,\theta_2,\theta_3, \theta_{\infty})$ is an
equivalent (Hamiltonian) form of ${\rm
P_{VI}}(\alpha,\beta,\gamma, \delta)$, where
$$
\alpha=\frac12\theta_{\infty}^2,\quad
\beta=-\frac12\theta_2^2,\quad \gamma=\frac12\theta_3^2,\quad
\delta=\frac12(1-\theta_1^2)
$$
(see \cite{Ok}).

There exist classical results (R.\,Fuchs \cite{Fu1}, R.\,Garnier
\cite{Gar}) on the connection of scalar Fuchsian equations of
second order with ${\rm P_{VI}}$ equations and Garnier systems.
Let us consider a scalar Fuchsian equation of second order with
singular points $a_1,\ldots,a_n$, $a_{n+1}=0$, $a_{n+2}=1$,
$a_{n+3}=\infty$ and apparent singularities $u_1,\ldots,u_n$ whose
Riemann scheme has the form
$$
\left(\begin{array}{ccc} a_i                 & \infty & u_k \\
                          0                  & \alpha &  0  \\
      \theta_i      & \alpha+\theta_{\infty} & 2
      \end{array}\right), \qquad i=1,\ldots,n+2,\;k=1,\ldots,n,\qquad\theta_i\not\in{\mathbb Z}
$$
($\alpha$ depends on the parameters $\theta_i$ according to the
classical Fuchs relation
$\sum_{i=1}^{n+2}\theta_i+\theta_{\infty}+2\alpha+2n=2n+1)$. There
is freedom of choice of such an equation. Its coefficients
$b_1(z)$, $b_2(z)$ depend on $a,u,\theta$ and $n$ arbitrary
parameters $v_1,\ldots,v_n$ ($v_i={\rm res}_{u_i}b_2(z)$).

Fix a set $\theta$ ($\theta_i\not\in{\mathbb Z}$) and consider an
($n$-dimensional) integral manifold $M$ of the system ${\cal
G}_n(\theta)$. Due to Theorem 4.1.2 from \cite{GP}, one has that
{\it Fuchsian equations corresponding to points $(a,u,v)\in M$
have the same monodromy\footnote{This property is defined
precisely in \S3.}. Inversely, points $(a,u,v)$ corresponding to
Fuchsian equations with the same mono\-dro\-my lie on the integral
manifold of the system ${\cal G}_n(\theta)$}.

Using the above relationship between Fuchsian and non-linear
differential equations one can deduce the known properties of the
latter as well as some new ones. In particular, we consider the
Riemann problem for some types of ${\rm GL}(2,\mathbb
C)$-representations (Theorems 3, 4) giving detailed proofs of the
statements from \cite{Go2}, and study movable singularities of
Garnier systems (Theorem 5).
\bigskip

{\bf \S2. Method of solution of the Riemann--Hilbert problem}
\medskip

In the study of problems related to the Riemann--Hilbert problem a
very useful tool is provided by linear gauge transformations of
the form
\begin{eqnarray}\label{transf}
y'=\Gamma(z)\,y
\end{eqnarray}
of the unknown function $y(z)$. The transformation (\ref{transf})
is said to be {\it holomorphically $($meromorphi\-cal\-ly$)$
invertible} at some point $z=a$, if the matrix $\Gamma(z)$ is
holomorphic (meromorphic) at this point and $\det\Gamma(a)\ne0$
($\det\Gamma(z)\not\equiv0$). This transformation transforms the
system (\ref{syst}) into the system
\begin{eqnarray}\label{transfmatr}
\frac{dy'}{dz}=B'(z)\,y', \qquad
B'(z)=\frac{d\Gamma}{dz}\Gamma^{-1}+
\Gamma\left(\sum_{i=1}^n\frac{B_i}{z-a_i}\right)\Gamma^{-1},
\end{eqnarray}
which is said to be, respectively, {\it holomorphically} or {\it
meromorphically equivalent} to the original system in a
neighbourhood of the point $a$.

An important property of meromorphic gauge transformations is the
fact that they do not change the monodromy (being meromorphic, the
matrix $\Gamma(z)$ is single-valued, therefore the ramification of
the fundamental matrix $\Gamma(z)Y(z)$ of the new system coincides
with the ramification of the matrix $Y(z)$).

Locally, in a neighbourhood of each point $a_k$, it is not
difficult to produce a system for which $a_k$ is a Fuchsian
singularity and the monodromy matrix at this point coincides with
the corresponding generator $G_k=\chi([\gamma_k])$ of the
representation (\ref{repr}). This system is
\begin{eqnarray}\label{systloc}
\frac{dy}{dz}=\frac{E_k}{z-a_k}\,y, \qquad E_k=\frac1{2\pi i}\ln
G_k,
\end{eqnarray}
with fundamental matrix $(z-a_k)^{E_k}:=e^{E_k\ln(z-a_k)}$. The
brunch of the logarithm of the matrix $G_k$ is chosen such that
the eigenvalues $\rho_k^{\alpha}$ of the matrix $E_k$ satisfy the
condition
\begin{eqnarray}\label{normcond}
0\leqslant{\rm Re}\,\rho_k^{\alpha}<1.
\end{eqnarray}
Indeed,
$$
\frac d{dz}(z-a_k)^{E_k}=\frac{E_k}{z-a_k}\,(z-a_k)^{E_k},
$$
and a single circuit around the point $a_k$ counterclockwise
transforms the matrix $(z-a_k)^{E_k}$ into the matrix
$$
e^{E_k(\ln(z-a_k)+2\pi i)}=e^{E_k\ln(z-a_k)}e^{2\pi iE_k}=
(z-a_k)^{E_k}G_k.
$$

Of course, not any system with the Fuchsian singularity $a_k$ and
the local monodromy matrix $G_k$ is holomorphically equivalent to
the system (\ref{systloc}) in a neighbourhood of this point.

Let $S_k$ be a non-singular matrix reducing the matrix $E_k$ to a
block-diagonal form $E'_k=S_kE_kS^{-1}_k={\rm diag}(E_k^1,\ldots,
E_k^m)$, where each block $E_k^j$ is an upper-triangular matrix
with the unique eigenvalue $\rho_k^j$. Consider a diagonal
integer-valued matrix $\Lambda_k={\rm diag}(\Lambda_k^1,\ldots,
\Lambda_k^m)$ with the same block structure and such that the
diagonal elements of each block $\Lambda_k^j$ form a
non-increasing sequence. Then according to (\ref{transfmatr}) the
transformation
$$
y'=\Gamma(z)\,y,\qquad \Gamma(z)=(z-a_k)^{\Lambda_k}S_k,
$$
transforms the system (\ref{systloc}) into the system
\begin{eqnarray}\label{systlocl}
\frac{dy'}{dz}=\left(\frac{\Lambda_k}{z-a_k}+
(z-a_k)^{\Lambda_k}\frac{E'_k}{z-a_k}(z-a_k)^{-\Lambda_k}\right)y',
\end{eqnarray}
for which the point $a_k$ is also a Fuchsian
singularity\footnote{As follows from the form of the matrices
$\Lambda_k$ and $E'_k$, the matrix
$(z-a_k)^{\Lambda_k}E'_k(z-a_k)^{-\Lambda_k}$ is holomorphic.} and
the matrix $G_k$ is the monodromy matrix. We call a set
$\{\Lambda_1,\ldots,\Lambda_n,S_1,\ldots,S_n\}$ of matrices having
the properties described above, a set of {\it admissible
matrices}.

According to Levelt's theorem \cite{Le}, the Fuchsian system
(\ref{syst}) is {\it holomorphically} equivalent to a system of
form (\ref{systlocl}) (with some matrix $\Lambda_k$) in a
neighbourhood of the singular point $a_k$, i.~e., the system has a
fundamental matrix
$$
Y_k(z)=U_k(z)(z-a_k)^{\Lambda_k}(z-a_k)^{E'_k},
$$
where the matrix $U_k(z)$ is holomorphically invertible at
$z=a_k$. The matrix $Y_k(z)$ is called the {\it Levelt fundamental
matrix} (its columns form the {\it Levelt basis}).

The eigenvalues $\beta_k^j$ of the residue matrix $B_k$ are said
to be the {\it exponents} of the Fuchsian system (\ref{syst}) at
the point $a_k$. They are invariants of the {\it holomorphic}
equivalence class of this system. From (\ref{systlocl}) it follows
that the exponents coincide with the eigenvalues of the matrix
$\Lambda_k+E'_k$. The matrix $\Lambda_k$ is said to be the {\it
valuation matrix} of the Fuchsian system (\ref{syst}) at the
singularity $a_k$. According to (\ref{normcond}), its diagonal
elements coincide with the integer parts of the numbers ${\rm
Re}\,\beta_k^j$.



The Riemann--Hilbert problem has a positive solution if one can
pass from the local systems (\ref{systlocl}) to a global Fuchsian
system defined on the whole Riemann sphere. The use of holomorphic
vector bundles and meromorphic connections proves to be effective
in the study of this question.

From the representation (\ref{repr}) one constructs over the
Riemann sphere a family ${\cal F}$ of holomorphic vector bundles
of rank $p$ with logarithmic (Fuchsian) connections having the
prescribed singular points $a_1,\ldots,a_n$ and monodromy
(\ref{repr}). The Riemann--Hilbert problem for the fixed
representation (\ref{repr}) is solved in the affirmative if some
bundle in the family ${\cal F}$ turns out to be holomorphically
trivial (then the corresponding logarithmic connection defines a
Fuchsian system with the given singularities $a_1,\ldots,a_n$ and
monodromy (\ref{repr}) on the whole Riemann sphere). We now
briefly present the construction of the family $\cal F$ (see
details in \cite{AB}, Sect. 3.1, 3.2 and 5.1).
\medskip

{\bf 1.} First, from the representation (\ref{repr}) over the
punctured Riemann sphere $B=\overline{\mathbb C}\setminus
\{a_1,\ldots,a_n\}$ one constructs a holomorphic vector bundle $F$
of rank $p$ with a holomorphic connection $\nabla$ that has the
given monodromy (\ref{repr}). The bundle $F$ over $B$ is obtained
from the holomorphically trivial bundle $\widetilde
B\times{\mathbb C}^p$ over the universal cover $\widetilde B$ of
the punctured Riemann sphere after identifications of the form
$(\tilde z, y)\sim(\sigma\tilde z, \chi(\sigma)y)$, where $\tilde
z\in\widetilde B$, $y\in{\mathbb C}^p$ and $\sigma$ is an element
of the group of deck transformations of $\widetilde B$ which is
identified with the fundamental group $\pi_1(B)$. Thus,
$F=\widetilde B\times{\mathbb C}^p/\sim$ and $\pi: F
\longrightarrow B$ is the natural projection. It is not difficult
to show that a gluing cocycle $\{g_{\alpha\beta}\}$ of the bundle
$F$ is defined by constant matrices $g_{\alpha\beta}$ after some
choice of a covering $\{U_{\alpha}\}$ of the punctured Riemann
sphere.

The holomorphic connection $\nabla$ can now be given by the set
$\{\omega_{\alpha}\}$ of matrix differential 1-forms
$\omega_{\alpha}\equiv0$, which obviously satisfy the gluing
conditions
\begin{equation}\label{glue}
\omega_{\alpha}=(dg_{\alpha\beta})g_{\alpha\beta}^{-1}+
g_{\alpha\beta}\omega_{\beta}g_{\alpha\beta}^{-1}
\end{equation}
on the intersections $U_{\alpha}\cap U_{\beta}\ne\varnothing$.
Furthermore, it follows from the construction of the bundle $F$
that the monodromy of the connection $\nabla$ coincides with
$\chi$.

{\bf 2.} Next, the pair $(F,\nabla)$ is extended to a bundle $F^0$
with a logarithmic connection $\nabla^0$ over the whole Riemann
sphere. For this, the set $\{U_{\alpha}\}$ should be supplemented
by small neighbourhoods $O_1,\ldots,O_n$ of the points
$a_1,\ldots,a_n$ respectively. An extension of the bundle $F$ to
each point $a_i$ looks as follows. For some non-empty intersection
$O_i\cap U_{\alpha}$ one takes $g_{i\alpha}(z)=(z-a_i)^{E_i}$ on
this intersection. For any other neighbourhood $U_{\beta}$ that
intersects $O_i$ one defines $g_{i\beta}(z)$ as the analytic
continuation of the matrix function $g_{i\alpha}(z)$ into $O_i\cap
U_{\beta}$ along a suitable path (so that the set
$\{g_{\alpha\beta}, g_{i\alpha}(z)\}$ defines a cocycle for the
covering $\{U_{\alpha}, O_i\}$ of the Riemann sphere). An
extension of the connection $\nabla$ to each point $a_i$ is given
by the matrix differential 1-form $\omega_i=E_idz/(z-a_i)$, which
has a simple pole at this point. Then the set $\{\omega_{\alpha},
\omega_i\}$ defines a logarithmic connection $\nabla^0$ in the
bundle $F^0$, since along with the conditions (\ref{glue}) for
non-empty $U_{\alpha}\cap U_{\beta}$, the conditions
$$
(dg_{i\alpha})g_{i\alpha}^{-1}+g_{i\alpha}\omega_{\alpha}g_{i\alpha}^{-1}
=\frac{E_i}{z-a_i}dz=\omega_i,\qquad O_i\cap
U_{\alpha}\ne\varnothing,
$$
also hold (see (\ref{systloc})). The pair $(F^0,\nabla^0)$ is
called the {\it canonical extension} of the pair $(F,\nabla)$.

{\bf 3.} In a way similar to that for the construction of the pair
$(F^0,\nabla^0)$, one can construct the family $\cal F$ of bundles
$F^{\Lambda}$ with logarithmic connections $\nabla^{\Lambda}$
having the given singularities $a_1,\ldots,a_n$ and monodromy
(\ref{repr}). For this, the matrices $g_{i\alpha}(z)$ in the
construction of the pair $(F^0,\nabla^0)$ should be replaced by
the matrices
$$
g^{\Lambda}_{i\alpha}(z)=(z-a_i)^{\Lambda_i}S_i(z-a_i)^{E_i},
$$
and the forms $\omega_i$ by the forms
$$
\omega^{\Lambda}_i=\left(\Lambda_i+(z-a_i)^{\Lambda_i}E'_i(z-a_i)^{-\Lambda_i}\right)
\frac{dz}{z-a_i},
$$
where $\{\Lambda_1,\ldots,\Lambda_n,S_1,\ldots,S_n\}$ are all
possible sets of admissible matrices. Then the conditions
\begin{eqnarray}\label{gluei}
(dg^{\Lambda}_{i\alpha})(g^{\Lambda}_{i\alpha})^{-1}+
g^{\Lambda}_{i\alpha}\omega_{\alpha}(g^{\Lambda}_{i\alpha})^{-1}
=\omega^{\Lambda}_i
\end{eqnarray}
again hold on the non-empty intersections $O_i\cap U_{\alpha}$
(see (\ref{systlocl})).
\medskip

{\bf Remark 1.} Strictly speaking, the bundle $F^{\Lambda}$ also
depends on the set $S=\{S_1,\ldots,S_n\}$ of the matrices $S_i$
reducing the monodromy matrices $G_i$ to an upper-triangular form.
In view of this dependence the bundles of the family $\cal F$
should be denoted by $F^{\Lambda, S}$. But in the following two
cases all bundles $F^{\Lambda, S}$ with a fixed $\Lambda$ are
holomorphically equivalent.

i) All points $a_i$ are {\it non-resonant}, i.~e., for each
valuation matrix $\Lambda_i$ all its blocks $\Lambda_i^j$ are
scalar matrices.

ii) Resonant points exist, but each resonant point $a_i$ has the
following property: for its monodromy matrix $G_i$ and any
$\lambda\in{\mathbb C}$ one has the inequality ${\rm
rank}(G_i-\lambda I)\geqslant p-1$.

In particular, in the two-dimensional case ($p=2$) if all
monodromy matrices $G_i$ are non-scalar, then bundles of the
family ${\cal F}$ depend on sets $\Lambda$ only.
\medskip

The exponents $\beta_i^j$ of the local Fuchsian system
$dy=\omega^{\Lambda}_iy$ are called the {\it exponents of the
logarithmic connection} $\nabla^{\Lambda}$ at the point $z=a_i$.

According to the Birkhoff--Grothendieck theorem, every holomorphic
vector bundle $E$ of rank $p$ over the Riemann sphere is
equivalent to a direct sum
$$
E\cong{\cal O}(k_1)\oplus\ldots\oplus{\cal O}(k_p)
$$
of line bundles which has a coordinate description of the form
$$
\left(U_0=\mathbb C, U_{\infty}=\overline{\mathbb
C}\setminus\{0\},g_{0\infty}=z^K\right), \qquad K={\rm
diag}\,(k_1,\ldots,k_p),
$$
where $k_1\geqslant\ldots\geqslant k_p$ is a set of integers which
is called the {\it splitting type} of the bundle $E$. The bundle
$E$ is holomorphically trivial if and only if it has the zero
splitting type.

The number $\deg E=\sum_{i=1}^pk_i$ equals the {\it degree} of the
bundle $E$. For the pair $(F^{\Lambda}, \nabla^{\Lambda})$ the
degree of the bundle $F^{\Lambda}$ coincides with the sum
$\sum_{i=1}^n\sum_{j=1}^p\beta_i^j=\sum_{i=1}^n{\rm tr}
(\Lambda_i+E_i)$ of the exponents of the connection
$\nabla^{\Lambda}$.

If some bundle $F^{\Lambda}$ in the family $\cal F$ is
holomorphically trivial then the corresponding logarithmic
connection $\nabla^{\Lambda}$ defines a global Fuchsian system
(\ref{syst}) that solves the Riemann--Hilbert problem. On the
other hand, in view of Levelt's theorem mentioned above, the
existence of a Fuchsian system with the given singular points
$a_1,\ldots,a_n$ and monodromy (\ref{repr}) implies the triviality
of some bundle in the family $\cal F$.

Thus, {\it the Riemann--Hilbert problem is soluble if and only if
at least one of the bundles of the family $\cal F$ is
holomorphically trivial} (see \cite{AB}, Th. 5.1.1).
\\

{\bf \S3. Isomonodromic deformations of Fuchsian systems}
\medskip

Let us include a Fuchsian system
\begin{eqnarray}\label{syst0}
\frac{dy}{dz}=\left(\sum_{i=1}^n\frac{B^0_i}{z-a^0_i}\right)y,
\qquad\sum_{i=1}^nB^0_i=0,
\end{eqnarray}
of $p$ equations into a family
\begin{eqnarray}\label{fam}
\frac{dy}{dz}=\left(\sum_{i=1}^n\frac{B_i(a)}{z-a_i}\right)y,
\qquad\sum_{i=1}^nB_i(a)=0,\qquad B_i(a^0)=B^0_i,
\end{eqnarray}
of Fuchsian systems holomorphically depending on the parameter
$a=(a_1,\ldots,a_n)\in D(a^0)$, where $D(a^0)$ is a disk of small
radius centered at the point $a^0=(a_1^0,\ldots,a_n^0)$ of the
space ${\mathbb C}^n\setminus\bigcup_{i\ne j}\{a_i=a_j\}$.

One says that the family (\ref{fam}) is {\it isomonodromic} (or it
is an {\it isomonodromic deformation} of the system
(\ref{syst0})), if for all $a\in D(a^0)$ the monodromies
$$
\chi: \pi_1(\overline{\mathbb C}\setminus\{a_1,\ldots,a_n\})
\longrightarrow{\rm GL}(p,{\mathbb C})
$$
of the corresponding systems are the same\footnote{Under small
variations of the parameter $a$ there exist canonical isomorphisms
of the fundamental groups $\pi_1(\overline{\mathbb
C}\setminus\{a_1,\ldots,a_n\})$ and $\pi_1(\overline{\mathbb
C}\setminus\{a_1^0,\ldots,a_n^0\})$ generating canonical
isomorphisms
$$
{\rm Hom}\left(\pi_1(\overline{\mathbb C}
\setminus\{a_1,\ldots,a_n\}), {\rm GL}(p,{\mathbb C})\right)/{\rm
GL}(p,{\mathbb C})\cong{\rm Hom}\left(\pi_1(\overline{\mathbb
C}\setminus\{a_1^0,\ldots,a_n^0\}),{\rm GL}(p,{\mathbb C})\right)/
{\rm GL}(p,{\mathbb C})
$$
of the spaces of conjugacy classes of representations for the
above fundamental groups; this allows one to compare $\chi$ for
various $a\in D(a^0)$.}. This means that for every value of $a$
from $D(a^0)$ there exists a fundamental matrix $Y(z, a)$ of the
corresponding system from (\ref{fam}) that has the same monodromy
matrices for all $a\in D(a^0)$. This matrix $Y(z, a)$ is called an
{\it isomonodromic fundamental matrix}.

For any isomonodromic family (\ref{fam}) there exists an
isomonodromic fundamental matrix that analytically depends on both
variables $z$ and $a$. An isomonodromic deformation preserves not
only the monodromy but also the exponents of the initial system
(thus, the eigenvalues of the residue matrices $B_i(a)$ of the
family (\ref{fam}) do not depend on the parameter $a$; see
\cite{Bo2} on the two latter statements).

Is it always possible to include the system (\ref{syst0}) into an
isomonodromic family of Fuchsian systems? The answer is positive.
Exactly, if the matrices $B_i(a)$ satisfy the {\it Schlesinger
equation} \cite{Sch}
$$
dB_i(a)=-\sum_{j=1, j\ne i}^n\frac{[B_i(a), B_j(a)]}
{a_i-a_j}\,d(a_i-a_j),
$$
then the family (\ref{fam}) is isomonodromic (in this case it is
called the {\it Schlesinger isomonodromic family}).

A Schlesinger isomonodromic family has the following property:
connection matrices between some fixed isomonodromic fundamental
matrix $Y(z,a)$ and local Levelt's bases at singular points do not
depend on $a$. Among all isomonodromic deformations of Fuchsian
systems with this property, the Schlesinger ones are distinguished
by the condition
$\left(d_aY(z,a)\right)Y^{-1}(z,a)|_{z=\infty}\equiv0$ (see
\cite{Bo2}).

It is well known that {\it for arbitrary initial conditions
$B_i(a^0)=B_i^0$ the Schlesinger equation has a unique solution
$\{B_1(a),\ldots,B_n(a)\}$ in some disk $D(a^0)$, and the matrices
$B_i(a)$ can be extended to the universal cover $Z$ of the space
${\mathbb C}^n\setminus \bigcup_{i\ne j}\{a_i=a_j\}$ as
meromorphic functions} (Malgrange's theorem \cite{Ma}). Thus, the
Schlesinger equation satisfies the Painlev\'e property.
\medskip

Here we should recall some facts on complex analytic sets and
meromorphic functions of several complex variables (one may see
details, for instance, in \cite{Ch}). We will need the former of
codimension one.

An analytic set $A\subset Z$ of codimension one (or complex
hypersurface) is defined locally as a set of zeros of a
holomorphic function, i.~e., near a point $a^0\in A$ it is defined
by an equation $f(a)=0$, where $f(a)$ is a holomorphic function in
a neighbourhood of the point $a^0$.

According to the Weierstrass preparation theorem, a holomorphic
function $f(a)$ vanishes at $a=a^0$ as a polynomial
$$
P(a)=(a_n-a_n^0)^k+c_1(a')(a_n-a_n^0)^{k-1}+\ldots+c_k(a')
$$
with coefficients $c_i(a')$ holomorphically depending on
$a'=(a_1,\ldots,a_{n-1})$ and vanishing at
$(a_1^0,\ldots,a_{n-1}^0)$, where $k\geqslant1$ is the order of
zero of $f(a_1^0,\ldots,a_{n-1}^0,a_n)$ at the point $a_n=a_n^0$
(we assume that $f(a_1^0,\ldots,a_{n-1}^0,a_n)\not\equiv0$). The
polynomial $P(a)$ is called the {\it Weierstrass polynomial} of
the function $f(a)$ at the point $a^0$.

The function $f(a)$ is called {\it irreducible} at the point
$a^0$, if it can not be decomposed into the product
$f(a)=f_1(a)f_2(a)$ of two functions holomorphic at this point and
vanishing there (it is equivalent to the irreducibility of its
Weierstrass polynomial). Any function $g(a)$, holomorphic at the
point $a^0$ and vanishing there, can be decomposed into the
product $g(a)=g_1^{m_1}(a)\ldots g_r^{m_r}(a)$, $m_i\in{\mathbb
N}$, of irreducible functions (up to holomorphic and not vanishing
at $a^0$ factors). One has $m_1=\ldots=m_r=1$, if and only if the
discriminant of the Weierstrass polynomial of the function $g(a)$
is not equal to zero identically.

If $f(a)$ can be chosen so that $df(a^0)\not\equiv0$, the point
$a^0$ is called a {\it regular point} of the set $A$, otherwise it
is called {\it critical}. If the function $f(a)$ is irreducible at
the point $a^0$, so is the set $A$ (i.~e., in any small
neighbourhood $U$ of this point $A\cap U$ can not be presented as
a union of analytic sets different from $A\cap U$); the inverse is
not always true. Thus, the set of points where $A$ is reducible,
is contained in the set of its critical points. In a neighbourhood
of any point the set $A$ can be presented as a union of
hypersurfaces irreducible at this point. The set $A$ is
irreducible at the point $a^0$ if and only if regular points of
$A\cap U$ form a connected set.
\medskip

{\bf Example 1.} a) The set $\{a\in{\mathbb C}^3\,|\,
a_1^2a_2^2-a_3^2=0\}=\{a_1a_2-a_3=0\}\bigcup\{a_1a_2+a_3=0\}$ is
reducible at any critical point (critical points of this set have
the form $(a_1,0,0)$ or $(0,a_2,0)$).
\smallskip

b) The set $A=\{a\in{\mathbb C}^3\,|\,a_1a_2-a_3^2=0\}$ is
irreducible at its critical point $a=0$, since the set
$A^0=(A\setminus\{0\})\bigcap\{|a_i|<\varepsilon^2\}$ is
connected. Indeed, $A$ is the image of ${\mathbb C}^2$ under the
transformation $(\xi,\eta)\mapsto (\xi^2,\eta^2,\xi\eta)$, and the
pre-image $\{|\xi|<\varepsilon, |\eta|<\varepsilon\}
\setminus\{0\}$ of the set $A^0$ is connected.
\smallskip

c) The set $A=\{a\in{\mathbb C}^3\,|\,a_1a_2^2-a_3^2=0\}$ is
reducible at the critical points $(a_1^0,0,0)$, $a_1^0\ne0$ (there
it can be presented as a union $\{a_2\sqrt{a_1}-a_3=0\}\bigcup
\{a_2\sqrt{a_1}+a_3=0\}$, where $\sqrt{a_1}$ is one of the
branches of the root near the point $a_1^0\ne0$) and irreducible
at the critical point $a=0$. Indeed, the set
$A^0=(A\setminus\{a_2=a_3=0\})\bigcap\{|a_i|<\varepsilon^2\}$ is
connected ($A$ is the image of ${\mathbb C}^2$ under the
transformation $(\xi,\eta)\mapsto (\xi^2,\eta,\xi\eta)$, and the
pre-image $\{|\xi|<\varepsilon, |\eta|<\varepsilon^2\}
\setminus\{\eta=0\}$ of the set $A^0$ is connected).
\medskip

A function $f(a)$ is {\it meromorphic} on $Z$, if it is
holomorphic on $Z\setminus P$, can not be extended to $P$
holomorphically and is presented as a quotient
$f(a)=\varphi(a)/\psi(a)$ of holomorphic functions in a
neighbourhood of every point $a^0\in P$ (hence, $\psi(a^0)=0$).
Thus, $P\subset Z$ is an analytic set of codimension one (it is
defined locally by the equation $\psi(a)=0$), which is called the
{\it polar locus} of the meromorphic function $f$. The points of
this set is divided into {\it poles} (at which the function
$\varphi$ does not vanish) and {\it ambiguous points} (at which
$\varphi=0$).

One can also define a divisor of a meromorphic function. Denote by
$A=N\cup P$ the union of the set $N$ of zeros and polar locus $P$
of the function $f$. Any regular point $a^0$ of the set $A$ can
belong to only one irreducible component of $N$ or $P$. Thus, one
can define the {\it order} of this component as the degree (taken
with "+", if $a^0\in N$, and with "$-$", if $a^0\in P$) of the
corresponding factor in the decomposition of the function
$\varphi$ or $\psi$ into irreducible factors. Then the {\it
divisor of the meromorphic function} $f$ is the pair $(A,
\kappa)$, where $\kappa=\kappa(a)$ is an integer-valued function
on the set of regular points of $A$ (which takes a constant value
on each its irreducible component, this value is equal to the
order of a component).
\medskip

{\bf Notation.} For the polar locus $P$ of the function $f$, and
$a^0\in P$, let us denote by $\Sigma_{a^0}(f)$ the sum of orders
of all irreducible components of $P\cap D(a^0)$.
\medskip

{\bf Example 2.} a) The function $f(a)=1/a_1a_2$ is meromorphic on
${\mathbb C}^2$. Its polar locus is $P=\{a_1a_2=0\}$ (all points
are poles), and the order of each component $\{a_i=0\}$ is equal
to $-1$, $\Sigma_0(f)=-2$.
\smallskip

b) The function $g(a)=a_1/a_2$ is meromorphic on ${\mathbb C}^2$.
Its polar locus is $P=\{a_2=0\}$ ($0$ is an ambiguous point, all
the others are poles), the set of zeros is $N=\{a_1=0\}\setminus
\{0\}$. The order of the component $\{a_1=0\}$ is equal to $1$,
the order of the component $\{a_2=0\}$ is equal to $-1$,
$\Sigma_0(g)=-1$.
\medskip

Let us return to the Schlesinger equation. The polar locus
$\Theta\subset Z$ of the extended matrix functions
$B_1(a),\ldots,B_n(a)$ is called the {\it Malgrange
$\Theta$-divisor}\footnote{In view of the above definition of a
divisor, here the term "divisor" is not precise enough.} ($\Theta$
depends on the initial conditions $B_i(a^0)=B_i^0$). Near a point
$a^*\in\Theta$ it is defined by the equation $\tau^*(a)=0$, where
$\tau^*(a)$ is a holomorphic function in a neighbourhood of the
point $a^*$ called a {\it local $\tau$-function} of the
Schlesinger equation. According to Miwa's theorem (\cite{JM}, see
also \cite{Bo3}) there exists a function $\tau(a)$ holomorphic on
the whole space $Z$ whose set of zeros coincides with $\Theta$. In
a neighbourhood of the point $a^*\in\Theta$ the global
$\tau$-function differs from the local one by a holomorphic
non-zero multiplier, and
$$
d\ln\tau(a)=\frac12\sum_{i=1}^n\sum_{j=1, j\ne i}^n\frac{{\rm tr}
(B_i(a)B_j(a))}{a_i-a_j}\,d(a_i-a_j).
$$

If we consider a system of the family (\ref{fam}) as an equation
for horizontal sections of the logarithmic connection
$\nabla^{\Lambda}$ (with singularities $a_1,\ldots,a_n$) in the
trivial bundle $F^{\Lambda}$, then the set $\Theta$ corresponds to
those points, where the extension of $F^{\Lambda}$ is not
holomorphically trivial.
\medskip

{\bf Example 3.} Consider a family
\begin{eqnarray}\label{constfam}
\frac{dy}{dz}=\left(\sum_{i=1}^n\frac{B_i}{z-a_i}\right)y,
\qquad\sum_{i=1}^nB_i=0,
\end{eqnarray}
of Fuchsian systems with {\it constant pairwise commuting}
matrices $B_i$.

This family is isomonodromic. Its isomonodromic fundamental matrix
$Y(z,a)$ has the form
$$
Y(z,a)=(z-a_1)^{B_1}\ldots(z-a_n)^{B_n},
$$
and the monodromy matrices $G_k=e^{2\pi iB_k}$ do not depend on
$a$.

The other explanation is the following. The matrices $B_i$
evidently satisfy the Schlesinger equation, so the family
(\ref{constfam}) is a Schlesinger isomonodromic family. It is
defined on the whole space ${\mathbb C}^n\setminus \bigcup_{i\neq
j}\{a_i=a_j\}$, therefore $\Theta=\varnothing$. Note that the
$\tau$-function satisfies the equation
$$
d\ln\tau(a)=\frac12\sum_{i=1}^n\sum_{j=1, j\ne
i}^n\frac{\alpha_{ij}}{a_i-a_j}\,d(a_i-a_j), \qquad
\alpha_{ij}={\rm tr}(B_iB_j),
$$
i. e., $\tau(a)=\prod_{i<j}(a_i-a_j)^{\alpha_{ij}}$ is a
holomorphic non-zero function on $Z$.
\medskip

In what follows we will use the theorem describing a general
solution of the Schlesinger equation near the $\Theta$-divisor in
the case $p=2$.
\medskip

{\bf Theorem 1} (Bolibrukh \cite{Bo5},\cite{Bo6}). {\it If the
monodromy of the two-dimensional family $(\ref{fam})$ is
irreducible, then $\Sigma_{a^*}(B_i)\geqslant2-n$ for any
$a^*\in\Theta$ $(i=1,\ldots,n)$.}
\medskip

Further we present a simplified proof of this theorem based on the
technique of the paper \cite{Bo3}, but first we recall this
technique in the proof of Proposition 1 below.




Consider an {\it irreducible} two-dimensional representation
$$
\chi:\pi_1(\overline{\mathbb C}\setminus
\{a_1,\ldots,a_n\})\longrightarrow{\rm GL}(2,{\mathbb C}),
$$
$a\in D(a^0)\subset{\mathbb C}^n\setminus \bigcup_{i\neq
j}\{a_i=a_j\}$, and the family $\cal F$ of holomorphic vector
bundles with logarithmic connections constructed by the
representation $\chi$.
\medskip

{\bf Proposition 1.} {\it Let
$(F^{\Lambda},\nabla^{\Lambda})\in{\cal F}$ and $\deg
F^{\Lambda}=0$. Then for all $a\in D(a^0)$, may be, with the
exception of an analytic subset of codimension one, the bundle
$F^{\Lambda}$ is holomorphically trivial} ({\it i.~e., for almost
all $a\in D(a^0)$ there exists a Fuchsian system with the given
singular points $a_1,\ldots,a_n$, monodromy $\chi$ and set
$\Lambda$ of valuation matrices}).

{\bf Proof.}  Choose an arbitrary point
$a^*=(a_1^*,\ldots,a_n^*)\in D(a^0)$. Suppose the corresponding
bundle $F^{\Lambda}$ is not holomorphically trivial:
$$
F^{\Lambda}\cong{\cal O}(-k)\oplus{\cal O}(k), \qquad k\geqslant1.
$$
Let us show that the set of points $a$, for which the
corresponding bundle $F^{\Lambda}$ is not holomorphically trivial,
is given by an equation $\tau^*(a)=0$ in a neighbourhood of the
point $a^*$, where $\tau^*(a)\not\equiv0$ is a holomorphic
function.

Consider an auxiliary system
$$
\frac{dy}{dz}=\left(\sum_{i=1}^n\frac{B_i^*}{z-a_i^*}\right)y
$$
with the monodromy $\chi$, valuation matrices
$\Lambda_1,\ldots,\Lambda_n$ at the points $a_1^*,\ldots,a_n^*$
respectively but also with the apparent Fuchsian singularity at
the infinity.

As follows from Bolibrukh's permutation lemma (Lemma 2 from
\cite{Bo4}), a fundamental matrix of the constructed system has
the form $Y(z)=U(z)z^K$ near the infinity, where
$$
U(z)=I+U_1\frac1z+U_2\frac1{z^2}+\ldots,\qquad K={\rm diag}(-k,k).
$$
Therefore, the residue matrix at the infinity is equal to $-K$,
and $\sum_{i=1}^n{B_i^*}=K$.

We need the following proposition which will be also used in the
further.
\medskip

{\bf Proposition 2} (Bolibrukh \cite{Bo4}). {\it Consider the
Fuchsian system (\ref{syst}) with the singularities
$a_1,\ldots,a_n$, apparent singularity $\infty$, monodromy
(\ref{repr}) and set $\Lambda=\{\Lambda_1,\ldots, \Lambda_n\}$ of
valuation matrices; furthermore $\sum_{i=1}^nB_i=K'={\rm diag}
(k_1,\ldots,k_p)$, where $k_1\leqslant\ldots\leqslant k_p$ are
integers.

The matrix $K'$ defines the splitting type of the bundle
$F^{\Lambda}$ if and only if the transformation $y'=z^{-K'}y$
transforms this system into the system that is holomorphic at the
infinity.}
\medskip

Due to this proposition the transformation $y'=z^{-K}y$ transforms
our auxiliary system into the system that is holomorphic at the
infinity, hence
$$
U(z)z^K=z^KV(z)
$$
for some matrix $V(z)$ holomorphically invertible at the infinity.
The latter relation implies that the upper-right element
$u_1^{12}$ of the matrix $U_1$ equals zero.

Using the theorem of existence and uniqueness for the Schlesinger
equation, include the constructed Fuchsian system into the
Schlesinger isomonodromic family
\begin{eqnarray}\label{addschl}
\frac{dy}{dz}=\left(\sum_{i=1}^n\frac{B_i(a)}{z-a_i}\right)y,\quad
B_i(a^*)=B_i^*, \quad \sum_{i=1}^nB_i(a)=K.
\end{eqnarray}
As shown in \cite{Bo3}, there exists an isomonodromic fundamental
matrix $Y(z, a)$ of this family of the form
\begin{eqnarray}\label{addschlfund}
Y(z,a)=U(z,a)z^K,\qquad
U(z,a)=I+U_1(a)\frac1z+U_2(a)\frac1{z^2}+\ldots,
\end{eqnarray}
at the infinity, $U(z,a^*)=U(z)$ and
\begin{eqnarray}\label{relation}
\frac{\partial U_1(a)}{\partial a_i}=-B_i(a), \qquad i=1,\ldots,n.
\end{eqnarray}

Since the monodromy $\chi$ is irreducible, among the upper-right
elements $b_i^{12}(a)$ of the corresponding matrices $B_i(a)$
there exists at least one that is not identically zero. Hence, in
view of (\ref{relation}), the similar element $u_1^{12}(a)$ of the
matrix $U_1(a)$ does not equal zero identically (while
$u_1^{12}(a^*)=u_1^{12}=0$).

Further, whereas
$$
\frac{dY(z,a)}{dz}Y^{-1}(z,a)=\sum_{i=1}^n\frac{B_i(a)}{z-a_i}=
\frac1z\sum_{i=1}^n\frac{B_i(a)}{1-\frac{a_i}z},
$$
from (\ref{addschlfund}) one gets the relation
\begin{eqnarray*}
-U_1(a)\frac1{z^2}+o(z^{-2})+\left(I+U_1(a)\frac1z+o(z^{-1})\right)\frac
Kz=\\=\biggl(\frac Kz+
\biggl(\sum_{i=1}^nB_i(a)a_i\biggr)\frac1{z^2}+o(z^{-2})\biggr)
\biggl(I+U_1(a)\frac1z+o(z^{-1})\biggr).
\end{eqnarray*}
Hence,
$$
-U_1(a)+[U_1(a),K]=\sum_{i=1}^nB_i(a)a_i.
$$
Therefore,
$$
(2k-1)u_1^{12}(a)=\sum_{i=1}^n{b_i^{12}(a)a_i}.
$$

Denote by $b_1(a)$ the sum $\sum_{i=1}^n{b_i^{12}(a)a_i}$. Then
$$
b_1(a)=(2k-1)u_1^{12}(a)\not\equiv0,\qquad b_1(a^*)=0.
$$
Consider the matrix
$$
\Gamma'_1(z,a)=\left(\begin{array}{cc} 1 & 0 \\
                \frac{1-2k}{b_1(a)}z & 1
                     \end{array}\right),
$$
holomorphically invertible (in $z$) off the infinity. One can
directly check that the matrix $U'(z,a)=\Gamma'_1U(z,a)$ has the
form
$$
U'(z,a)=\left(U'_0(a)+U'_1(a)\frac1z+\ldots\right)z^{{\rm diag}
(1,-1)}, \qquad U'_0(a)=\left(\begin{array}{cc} 0 & \frac{b_1(a)}{2k-1} \\
                              \frac{1-2k}{b_1(a)} & \frac{f(a)}{b_1(a)}
                              \end{array}\right),
$$
where $f(a)$ is a holomorphic function at the point $a^*$. Thus,
the gauge transformation $y_1=\Gamma_1(z,a)y$,
$\Gamma_1(z,a)=U'_0(a)^{-1}\Gamma'_1(z,a)$, transforms a system of
the family (\ref{addschl}) into the Fuchsian system with the
fundamental matrix $Y^1(z,a)=\Gamma_1(z,a)Y(z,a)$ of the form
(\ref{addschlfund}) at the infinity (and does not change
valuations at the points $a_1,\ldots,a_n$), where all involved
matrices are equipped with the upper index 1, and $K^1={\rm
diag}(-k+1,k-1)$. This expansion is valid only in the exterior of
some analytic subset of codimension one which is the set of zeros
of the function $b_1(a)\not\equiv0$.

Note also that the transformed family is a Schlesinger
isomonodromic family. Indeed, its connection matrices do not
depend on $a$ (the transformation does not change those of the
Schlesinger family (\ref{addschl})), and
$(d_aY^1(z,a))Y^1(z,a)^{-1}|_{z=\infty}=(d_aU^1(z,a))U^1(z,a)^{-1}|_{z=\infty}
\equiv0$ according to the form of the matrix $U^1(z,a)$.

After $k$ steps of the above procedure of Bolibrukh we will get a
Fuchsian family holomorphic at the infinity. It is defined in a
neighbourhood of the point $a^*$ outside of the analytic subset
$\{\tau^*(a)=0\}$, $\tau^*(a)=b_1(a)\ldots b_k(a)$, where $b_j(a)$
appears at the $j$-th step of the Bolibrukh procedure in the same
way as $b_1(a)$ does. This means that for all
$a\not\in\{\tau^*(a)=0\}$ from the neighbourhood of the point
$a^*$ there exists a Fuchsian system with the singularities
$a_1,\ldots,a_n$, monodromy $\chi$ and set $\Lambda$ of valuation
matrices. $\hfill{\Box}$
\medskip

{\bf Definition.} Recall that if all generators $G_i$ of the {\it
two-dimensional} representation $\chi$ are non-scalar matrices,
then bundles of the family $\cal F$ depend on sets $\Lambda$ only
(see Remark 1). One calls such representations {\it non-smaller}.
In the opposite case, if $l$ monodromy matrices are scalar, $\chi$
is called $l$-{\it smaller}.
\medskip

{\bf Corollary 1.} {\it If $\chi$ is an irreducible non-smaller
${\rm SL}(2, {\mathbb C})$-representation with generators
$G_1,\ldots,G_n$, then for almost all $a\in D(a^0)$ there exists a
family} ({\it depending on the parameter ${\bf
m}=(m_1,\ldots,m_n)\in{\mathbb Z}_+^n$})
$$
\frac{dy}{dz}=\left(\sum_{i=1}^n\frac{B_i^{\bf
m}(a)}{z-a_i}\right)y, \qquad\sum_{i=1}^nB_i^{\bf m}(a)=0,
$$
{\it of Fuchsian systems with the singularities $a_1,\ldots,a_n$,
monodromy $\chi$ and exponents $\pm(m_k+\rho_k)$, where $\rho_k$
is one of the eigenvalues of the matrix $E_k=(1/2\pi i)\ln G_k$
$(k=1,\ldots,n)$. Furthermore, $B_n^{\bf m}(a)={\rm diag}
(m_n+\rho_n,-m_n-\rho_n)$ are diagonal matrices.}
\smallskip

{\bf Proof.} If a set $\Lambda=\{\Lambda_1,\ldots,\Lambda_n\}$ of
admissible matrices satisfies the conditions ${\rm
tr}(\Lambda_k+E_k)=0$, $k=1,\ldots,n$, then by Proposition 1 for
all $a\in D(a^0)$, may be, with the exception of an analytic
subset $\Theta_{\Lambda}$ of codimension one, the corresponding
bundle $F^{\Lambda}$ is holomorphically trivial and the
logarithmic connection $\nabla^{\Lambda}$ defines a Fuchsian
system with the singularities $a_1,\ldots,a_n$, monodromy $\chi$
and set $\Lambda$ of valuation matrices.

By the relations $e^{2\pi iE_k}=G_k$, $\det G_k=1$, the sum
$\rho_k^1+\rho_k^2$ of the eigenvalues of the matrix $E_k$ is an
integer, and it equals $0$ or $1$ by the condition
(\ref{normcond}). Fix an order of the eigenvalues $\rho_k^1$,
$\rho_k^2$ and put $\rho_k=\rho_k^1$.

1) If $\rho_k^1+\rho_k^2=0$, then one can take $\Lambda_k={\rm
diag} (m_k,-m_k)$, $m_k\in{\mathbb Z}_+$ (but if $\rho_n=0$, then
$m_n\in{\mathbb N}$).

2) If $\rho_k^1+\rho_k^2=1$, then one can take $\Lambda_k={\rm
diag} (m_k,-m_k-1)$, $m_k\in{\mathbb Z}_+$.

Thus, for all $a\in D(a^0)\setminus\Theta_{\bf m}$ the
representation $\chi$ can be realized by a Fuchsian system with
the singular points $a_1,\ldots,a_n$ and exponents
$\pm(m_1+\rho_1),\ldots, \pm(m_n+\rho_n)$. Moreover, the residue
matrix at the point $a_n$ is diagonalisable (because its
eigenvalues $\pm(m_n+\rho_n)$ do not equal zero by the
construction). Then the statement of the corollary is valid for
all $a\in D(a^0)\setminus\bigcup_{\bf m}\Theta_{\bf m}$.
\hfill{$\Box$}
\medskip

{\bf Proof of Theorem 1.} For $a^*\in\Theta$ the corresponding
vector bundle $F^{\Lambda}\cong{\cal O}(-k)\oplus{\cal O}(k)$ is
not holomorphically trivial and, as shown in the proof of
Proposition 1, the $\Theta$-divisor of the family (\ref{fam}) in a
neighbourhood of the point $a^*$ is the set of zeros of the
function $\tau^*(a)=b_1(a)\ldots b_k(a)$ constructed by the
auxiliary family (\ref{addschl}). Let us denote by $B_i^*(a)$ the
residue matrices of the latter (to tell them from those $B_i(a)$
of the initial family (\ref{fam})). They are holomorphic in a
neighbourhood of the point $a^*$.

The functions $b_j(a)$ are irreducible at $a^*$, since
$db_j(a^*)\not\equiv0$. For instance,
$$
db_1(a)=(2k-1)du_1^{12}(a)=(1-2k)\sum_{i=1}^nb_i^{12}(a)da_i
$$
in view of (\ref{relation}), and the equality $db_1(a^*)\equiv0$
implies $b_1^{12}(a^*)=\ldots=b_n^{12}(a^*)=0$, which contradicts
the irreducibility of the monodromy.

One can assume that $\tau^*(a)=b_1^{m_1}(a)\ldots b_r^{m_r}(a)$,
$m_1+\ldots+m_r=k$ (some factors are equal). Now let us show that
the order of each component $\{b_j(a)=0\}$ is not less than
$-2m_j$. It is sufficient to consider the first step of the
Bolibrukh procedure. The transformation $y_1=\Gamma_1(z,a)y$
transforms the auxiliary family into the family with the
coefficient matrix of the form
$$
\frac{d\Gamma_1}{dz}\Gamma_1^{-1}+\Gamma_1\left(
\sum_{i=1}^n\frac{B_i^*(a)}{z-a_i}\right)\Gamma_1^{-1},
$$
where
$$
\Gamma_1(z,a)=U'_0(a)^{-1}\Gamma'_1(z,a)=
\left(\begin{array}{cc} \frac{f(a)}{b_1(a)} & \frac{b_1(a)}{1-2k} \\
                        \frac{2k-1}{b_1(a)} & 0 \end{array}\right)
\left(\begin{array}{cc} 1 & 0 \\
     \frac{1-2k}{b_1(a)}z & 1
      \end{array}\right)=
\left(\begin{array}{cc} \frac{f(a)}{b_1(a)}+z & \frac{b_1(a)}{1-2k} \\
                        \frac{2k-1}{b_1(a)}   & 0
      \end{array}\right).
$$
Thus, the residue matrices $B_i^1(a)$ of the transformed family
have the form
$$
B_i^1(a)=\left(\begin{array}{cc} \frac{f(a)}{b_1(a)}+a_i & \frac{b_1(a)}{1-2k} \\
                        \frac{2k-1}{b_1(a)}   & 0
               \end{array}\right) B_i^*(a)
\left(\begin{array}{cc} \frac{f(a)}{b_1(a)}+a_i & \frac{b_1(a)}{1-2k} \\
                        \frac{2k-1}{b_1(a)}   & 0
      \end{array}\right)^{-1},
$$
i. e., the matrices $b_1(a)^2 B_i^1(a)$ are holomorphic in
$D(a^*)$.

After the final ($k$-th) step of the procedure we get the
Schlesinger isomonodromic family with the residue matrices
$B_i^k(a)$ which are simultaneously conjugated to the
corresponding $B_i(a)$ of the initial family (\ref{fam}) by some
constant matrix $S$ (this follows from the uniqueness of a
solution to the Schlesinger equation). Therefore,
$\Sigma_{a^*}(B_i)\geqslant -2m_1-\ldots-2m_r=-2k\geqslant2-n$
(see (\ref{ineq})). \hfill{$\Box$}
\medskip

In the case of dimension $p>2$ one can also apply a similar
procedure to find a local $\tau$-function $\tau^*(a)=b_1(a)\ldots
b_s(a)$. We can not assert that the functions $b_j(a)$ are
irreducible at the point $a^*$. But if for each $b_j(a)$ all its
irreducible factors are distinct (this is the case when the
discriminant of the Weierstrass polynomial of each $b_j(a)$ is not
identically zero), then one can estimate the order $\kappa$ of
each irreducible component of the $\Theta$-divisor as follows (see
\cite{Go}):
$$
\kappa\geqslant-\frac{(n-2)p(p-1)}2,
$$
if the monodromy of the family is irreducible, and
$$
\kappa\geqslant-\sum_{i=1}^n(M_i-\mu_i)\frac{p(p-1)}2
$$
in the case of reducible monodromy, where $\mu_i<M_i$ are integers
that bound real parts of the eigenvalues of the residue matrix
$B_i(a)$.

The following auxiliary lemma is a simplified version of
Proposition 6.4.1 from \cite{GP}.
\medskip

{\bf Lemma 1.} {\it Consider a two-dimensional Schlesinger
isomonodromic family of the form
$$
\frac{dy}{dz}=\left(\sum_{i=1}^n\frac{B_i(a)}{z-a_i}\right)y,\qquad
\sum_{i=1}^nB_i(a)=K={\rm diag}(\theta,-\theta),\quad
\theta\in{\mathbb C},
$$
and the function $b(a)=\sum_{i=1}^nb_i^{12}(a)a_i$, where
$b_i^{12}(a)$ are the upper-right elements of the matrices
$B_i(a)$ respectively. Then the differential of the function
$b(a)$ is given by the formula}
$$
db(a)=(2\theta+1)\sum_{i=1}^nb_i^{12}(a)da_i.
$$

Note that we can not directly apply calculations of Proposition 1,
because an isomonodromic fundamental matrix of the family not
necessary has the form (\ref{addschlfund}) (the monodromy at the
infinity can be non-diagonal).

{\bf Proof.} The differential $db(a)$ has the form
$$
db(a)=\sum_{i=1}^na_idb_i^{12}(a)+\sum_{i=1}^nb_i^{12}(a)da_i.
$$
To find the first of the two latter summands, let us use the
Schlesinger equation
$$
dB_i(a)=-\sum_{j=1, j\ne i}^n
\frac{[B_i(a),B_j(a)]}{a_i-a_j}\,d(a_i-a_j)
$$
for the matrices $B_i(a)$. Then we have
\begin{eqnarray*}
\sum_{i=1}^n a_i\,dB_i(a)&=&-\sum_{i=1}^n\sum_{j=1, j\ne i}^n
a_i\frac{[B_i(a),B_j(a)]}{a_i-a_j}\,d(a_i-a_j)=
-\sum_{i=1}^n\sum_{j>i}^n\,[B_i(a),B_j(a)]d(a_i-a_j)=\\
&=&-\sum_{i=1}^n\Bigl[B_i(a),\sum_{j=1, j\ne i}^n B_j(a)\Bigr]da_i
=-\sum_{i=1}^n[B_i(a),K]da_i.
\end{eqnarray*}

The upper-right element of the latter matrix 1-form is equal to
$\sum_{i=1}^n2\theta b_i^{12}(a)da_i$, hence
$\sum_{i=1}^na_idb_i^{12}(a)=2\theta\sum_{i=1}^n b_i^{12}(a)da_i$,
and $db(a)=(2\theta+1)\sum_{i=1}^nb_i^{12}(a)da_i$. \hfill{$\Box$}
\\

{\bf \S4. The Riemann--Hilbert problem and the Painlev\'e VI
equation}
\medskip

As mentioned earlier, the problem of constructing a Fuchsian
differential equation (\ref{urav}) with the given singularities
$a_1,\ldots,a_n$ and monodromy (\ref{repr}) has a negative
solution in general case. In the construction there arise apparent
singular points. In the case of {\it irreducible} representation
A.\,A.\,Bolibrukh \cite{Bo4} obtained the formula for the minimal
number of such singularities. It is given below.

We consider the family $\cal F$ of holomorphic vector bundles
$F^{\Lambda}$ with logarithmic connections $\nabla^{\Lambda}$
constructed from the representation (\ref{repr}). The {\it
Fuchsian weight} of the bundle $F^{\Lambda}$ is defined as the
quantity
$$
\gamma(F^{\Lambda})=\sum_{i=1}^p(k_1-k_i),
$$
where $(k_1,\ldots,k_p)$ is the splitting type of $F^{\Lambda}$.

If the representation (\ref{repr}) is {\it irreducible}, then the
splitting type of the bundle $F^{\Lambda}$ satisfies the
inequalities
\begin{eqnarray}\label{ineq}
k_i-k_{i+1}\leqslant n-2, \qquad i=1,\ldots,p-1
\end{eqnarray}
(see \cite{Bo4}, Cor. 3). Therefore, the quantity
$$
\gamma_{\max}(\chi)=\max_{F^{\Lambda}\in\cal
F}\,\gamma(F^{\Lambda})\leqslant\frac{(n-2)p(p-1)}2
$$
is defined for such a representation, and is called the {\it
maximal Fuchsian weight} of the irreducible representation $\chi$.

The minimal possible number $m_0$ of apparent singular points
emerging in the construction of a Fuchsian equation (\ref{urav})
with the {\it irreducible} monodromy (\ref{repr}), is given by the
formula
\begin{equation}\label{numbersing}
m_0=\frac{(n-2)p(p-1)}2-\gamma_{\max}(\chi).
\end{equation}
In the case of reducible representation there exists the estimate
$m_0\leqslant 1+(n+1)p(p-1)/2$ obtained in \cite{VG}.

In particular, it follows from the formula (\ref{numbersing}) that
{\it a set of singular points $a_1$, $a_2$, $a_3$ $(n=3)$ and
irreducible two-dimensional representation $(p=2)$ can always be
realized by a Fuchsian differential equation of second order},
since in this case $\gamma(F^{\Lambda})=1$ for any bundle
$F^{\Lambda}$ of odd degree.

A ${\rm P_{VI}}$ equation appears when one solves the problem of
constructing a Fuchsian differential equation of {\it second}
order with {\it four} given singularities and an {\it irreducible}
monodromy. Further we recall this fact.

Let us consider the four points $t, 0, 1, \infty$ ($t\in D(t^*)$,
where $D(t^*)\subset{\mathbb C}\setminus\{0, 1\}$ is a disk of
small radius centered at the point $t^*$) and an {\it irreducible
non-smaller} representation
\begin{eqnarray}\label{repr1}
\chi^*: \pi_1({\mathbb C}\setminus\{t,0,1\})\longrightarrow {\rm
GL}(2,\mathbb C)
\end{eqnarray}
generated by matrices $G_1, G_2, G_3$ corresponding to the points
$t, 0, 1$ (recall that in this case bundles of the family $\cal F$
depend on sets $\Lambda$ of valuation matrices only).

Depending on the location of the point $t$, there are two possible
cases.

{\it 1) Every vector bundle $F^{\Lambda}$ in the family $\cal F$
constructed with respect to the given four points and
representation $\chi^*$, such that $\deg F^{\Lambda}=0$, is
holomorphically trivial} (as follows from Proposition 1, this is
the case for almost all values $t\in D(t^*)$).

{\it 2) Among the elements of the family $\cal F$ there exists a
non-trivial holomorphic vector bundle $F^{\Lambda}$ of degree
zero.} (Denote by $\widetilde\Theta$ the set of values of the
parameter $t$ that correspond to this case.)

It follows from the inequalities (\ref{ineq}) that
$\gamma_{\max}(\chi^*)\leqslant2$, therefore in the first case the
splitting type of a non-trivial holomorphic vector bundle
$F^{\Lambda}$ (of non-zero degree) can be $(k, k-1)$ or $(k, k)$
only. The case $(k+1, k-1)$ is impossible, since then the bundle
$F^{\Lambda}\otimes{\cal O}(-k)$ constructed with respect to the
set of valuation matrices ${\Lambda_1-kI,\Lambda_2,\Lambda_3,
\Lambda_{\infty}}$ has degree zero, i.~e., is holomorphically
trivial, but at the same time its splitting type is $(1,-1)$.
Consequently, $\gamma_{\max}(\chi^*)=1$ in the first case.

In the second case the splitting type of the non-trivial
holomorphic vector bundle of degree zero equals $(1,-1)$, and
$\gamma_{\max}(\chi^*)=2$ in this case.

Thus, in view of the formula (\ref{numbersing}), for almost all
values $t\in D(t^*)$ the set of points $t, 0, 1, \infty$ and
representation $\chi^*$ can be realized by a Fuchsian differential
equation of second order with one apparent singularity. We denote
this singularity by $u(t)$ regarding it as a function of the
parameter $t$. It turns out that the function $u(t)$ satisfies the
equation (\ref{PVI}) for some values of the constants $\alpha,
\beta, \gamma, \delta$, if $\chi^*$ is an ${\rm SL}(2,{\mathbb
C})$-representation\footnote{The ${\rm P_{VI}}$ equation was
obtained by R.\,Fuchs precisely as a differential equation that is
satisfied by the apparent (fifth) singularity $\lambda(t)$ of some
Fuchsian equation of second order with the singular points $0, 1,
t, \infty$ and ${\rm SL}(2,{\mathbb C})$-monodromy independent of
the parameter $t$.}. Let us explain this interesting fact by using
isomonodromic deformations of Fuchsian systems.

By Corollary 1 we can choose a value $t=t^0\in D(t^*)$ for which
the representation $\chi^*$ is realized by Fuchsian systems
\begin{eqnarray}\label{fuchs2}
\frac{dy}{dz}=\left(\frac{B_1^{\bf m}}{z-t^0}+\frac{B_2^{\bf m}}
z+\frac{B_3^{\bf m}}{z-1}\right)y,\qquad {\bf m}=
(m_1,m_2,m_3,m_{\infty})\in{\mathbb Z}_+^4,
\end{eqnarray}
with the singular points $t^0,0,1,\infty$ (the eigenvalues of the
matrices $B_k^{\bf m}$ are $\pm(m_k+\rho_k)$, and the matrices
$B_{\infty}^{\bf m}=-B_1^{\bf m}-B_2^{\bf m}-B_3^{\bf m}$ are
diagonal).

Any system of the form (\ref{fuchs2}) can be included into the
Schlesinger isomonodromic family
\begin{eqnarray}\label{schl2}
\frac{dy}{dz}=\left(\frac{B_1^{\bf m}(t)}{z-t}+\frac{B_2^{\bf
m}(t)}z+\frac{B_3^{\bf m}(t)}{z-1}\right)y, \qquad B_k^{\bf
m}(t^0)=B_k^{\bf m},
\end{eqnarray}
of Fuchsian systems with the singularities $t,0,1,\infty$ which
depends holomorphically on the parameter $t\in D(t^0)$.
Furthermore, $B_1^{\bf m}(t)+B_2^{\bf m}(t)+B_3^{\bf
m}(t)=-B_{\infty}^{\bf m}={\rm diag}(-m_{\infty}-\rho_{\infty},
m_{\infty}+\rho_{\infty})$.

Denote by $B_{\bf m}(z,t)=(b_{ij}^{\bf m}(z,t))$ the coefficient
matrix of the family (\ref{schl2}). Since the upper-right element
of the matrix $B_1^{\bf m}(t)+B_2^{\bf m}(t)+B_3^{\bf
m}(t)=-B_{\infty}^{\bf m}$ is equal to zero, for every fixed $t$
the same element of the matrix $z(z-1)(z-t)B_{\bf m}(z,t)$ is a
polynomial of first degree in $z$. We define $\tilde u_{\bf m}(t)$
as the unique root of this polynomial. Next we use the following
theorem (see \cite{JM} or \cite{GP}, Cor. 6.2.2).
\medskip

{\bf Theorem 2.}  {\it The function $\tilde u_{\bf m}(t)$
satisfies the equation $(\ref{PVI})$, where the constants
$\alpha$, $\beta$, $\gamma$, $\delta$ are connected with the
parameter ${\bf m}=(m_1,m_2,m_3,m_{\infty})$ by the relations
$$
\alpha=\frac{(2m_{\infty}+2\rho_{\infty}-1)^2}2,\quad
\beta=-2(m_2+\rho_2)^2,\quad \gamma=2(m_3+\rho_3)^2,\quad
\delta=\frac12-2(m_1+\rho_1)^2.
$$
}

Let us consider the row vectors
$$
h_0^{\bf m}=(1,0),\qquad h_1^{\bf m}(z,t)=\frac{dh_0^{\bf m}}{dz}
+h_0^{\bf m}B_{\bf m}(z,t)=(b_{11}^{\bf m},b_{12}^{\bf m})
$$
and the matrix composed from them,
$$
\Gamma_{\bf m}(z,t)=\left(\begin{array}{c} h_0^{\bf m}\\
                                           h_1^{\bf m}
                          \end{array}\right)=
                    \left(\begin{array}{cc}
                           1 & 0 \\
                           b_{11}^{\bf m} & b_{12}^{\bf m}
                          \end{array}\right),
$$
which is meromorphically invertible on $\overline{\mathbb C}\times
D(t^0)$, since $\det\Gamma_{\bf m}(z,t)=b_{12}^{\bf
m}(z,t)\not\equiv0$ by the irreducibility of the representation
$\chi^*$. We define functions $p_{\bf m}(z,t)$ and $q_{\bf
m}(z,t)$, meromorphic on $\overline{\mathbb C}\times D(t^0)$, so
that the relation
$$
h_2^{\bf m}(z,t):=\frac{dh_1^{\bf m}}{dz}+h_1^{\bf m}B_{\bf m}
(z,t)=(-q_{\bf m},-p_{\bf m})\Gamma_{\bf m}(z,t)
$$
holds. Then
$$
\frac{d\Gamma_{\bf m}}{dz}=\frac d{dz}\left(\begin{array}{c}
                                  h_0^{\bf m} \\ 
                                  h_1^{\bf m}
                                 \end{array}\right)=
                           \left(\begin{array}{c}
                                  h_1^{\bf m} \\
                                  h_2^{\bf m}
                                 \end{array}\right)-
                           \left(\begin{array}{c}
                                  h_0^{\bf m} \\ h_1^{\bf m}
                                 \end{array}\right)B_{\bf m}(z,t)=
                           \left(\begin{array}{cc}
                                  0 & 1 \\
                                 -q_{\bf m} & -p_{\bf m}
                                 \end{array}\right)\Gamma_{\bf m}-
                                 \Gamma_{\bf m}B_{\bf m}(z,t),
$$
whence,
$$
\left(\begin{array}{cc}
        0 & 1 \\
       -q_{\bf m} & -p_{\bf m}
      \end{array}\right)=\frac{d\Gamma_{\bf m}}{dz}\Gamma_{\bf m}^{-1}
      +\Gamma_{\bf m}B_{\bf m}(z,t)\Gamma_{\bf m}^{-1}.
$$

The latter means that for every fixed $t\in D(t^0)$ the gauge
transformation $y'=\Gamma_{\bf m}(z,t)y$ transforms the
corresponding system of the family (\ref{schl2}) into the system
$$
\frac{dy'}{dz}=\left(\begin{array}{cc}
                       0 & 1 \\
                      -q_{\bf m} & -p_{\bf m}
                     \end{array}\right)y',
$$
the first coordinate of whose solution is the solution of the
scalar equation
\begin{eqnarray}\label{urav2}
\frac{d^2w}{dz^2}+p_{\bf m}(z,t)\frac{dw}{dz}+q_{\bf m}(z,t)w=0.
\end{eqnarray}
This (Fuchsian) equation has the singular points $t,0,1,\infty$
and monodromy $\chi^*$, but it also has the apparent singularity
$u_{\bf m}(t)$ which is a zero of the function $\det\Gamma_{\bf
m}(z,t)=b^{\bf m}_{12}(z,t)$, as follows from the construction of
the functions $p_{\bf m}(z,t)$, $q_{\bf m}(z,t)$. By Theorem 2 the
function $u_{\bf m}(t)$ satisfies an equation ${\rm P_{VI}}$.

Thus, we can formulate the following statement.
\medskip

{\bf Theorem 3.}

{\it $\rm i)$ The set of the points $t,0,1,\infty$ and any
irreducible non-smaller ${\rm SL}(2,{\mathbb C})$-representation
$(\ref{repr1})$ can be realized by the family $($depending on the
parameter ${\bf m}\in{\mathbb Z}_+^4)$ of scalar Fuchsian
equations $(\ref{urav2})$ with one apparent singularity\footnote{
The apparent singular point $u_{\bf m}(t)$ of every equation from
this family, as a function of the parameter $t\in D(t^*)$,
satisfies the equation ${\rm P_{VI}}$ with the constants $\alpha$,
$\beta$, $\gamma$, $\delta$ given by Theorem 2.}.

$\rm ii)$ The set $\widetilde\Theta\supset\bigcup_{\bf m}\{t\in
D(t^*)|u_{\bf m}(t)=t,0,1, or\; \infty\}$ is a countable set of
parameter values for which the Riemann--Hilbert problem for scalar
Fuchsian equations under consideration is soluble without apparent
singularities.}
\medskip

Being solutions of ${\rm P_{VI}}$ equations, the functions $u_{\bf
m}(t)$ have only poles as movable singularities (in other words,
they can be extended to the universal covering $H$ of the space
${\mathbb C}\setminus\{0,1\}$ as meromorphic functions). What one
can say about their pole orders?

Denote by $b_1^{\bf m}(t)$, $b_2^{\bf m}(t)$, $b_3^{\bf m}(t)$ the
upper-right elements of the matrices $B_1^{\bf m}(t)$, $B_2^{\bf
m}(t)$, $B_3^{\bf m}(t)$ respectively (recall that $b_1^{\bf
m}(t)+b_2^{\bf m}(t)+b_3^{\bf m}(t)\equiv0$). Since
$$
b^{\bf m}_{12}(z,t)=\frac{(tb_1^{\bf m}+b_3^{\bf m})z+tb_2^{\bf
m}}{z(z-1)(z-t)},
$$
the function $u_{\bf m}(t)$ is given by the relation
$$
(tb_1^{\bf m}+b_3^{\bf m})u_{\bf m}=-tb_2^{\bf m},
$$
from which it follows that poles of the function $u_{\bf m}(t)$
are poles of the function $b_2^{\bf m}(t)$ or zeros of the
function $tb_1^{\bf m}(t)+b_3^{\bf m}(t)$.

By Theorem 1 (with $n=4$) a pole order of the function $b_i^{\bf
m}(t)$ does not exceed two. Applying Lemma 1 to the family
(\ref{schl2}), where $(a_1,a_2,a_3)=(t,0,1)$, one gets
$$
\frac d{dt}(tb_1^{\bf m}(t)+b_3^{\bf m}(t))=
(-2m_{\infty}-2\rho_{\infty}+1)b_1^{\bf m}(t).
$$

If $(m_{\infty},\rho_{\infty})\ne(0,1/2)$, then
$\theta=-2m_{\infty}-2\rho_{\infty}+1\ne0$. In this case a pole of
the function $b_2^{\bf m}(t)$ is also a pole for $tb_1^{\bf
m}(t)+b_3^{\bf m}(t)$, since
$$
tb_1^{\bf m}(t)+b_3^{\bf m}(t)=b_1^{\bf m}(t)+b_3^{\bf
m}(t)+(t-1)b_1^{\bf m}(t)=-b_2^{\bf m}(t)+\frac{t-1}{\theta}\frac
d{dt}(tb_1^{\bf m}(t)+b_3^{\bf m}(t)).
$$
From this relation it also follows that any zero $t_0$ of the
function $tb_1^{\bf m}(t)+b_3^{\bf m}(t)$ can be simple only.
Indeed, if $t_0b_1^{\bf m}(t_0)+b_3^{\bf m}(t_0)=0$ and $\frac
d{dt}(tb_1^{\bf m}(t)+b_3^{\bf m}(t))|_{t=t_0}=0$, then $b_2^{\bf
m}(t_0)=0$ and $b_1^{\bf m}(t_0)=b_3^{\bf m}(t_0)=0$, which
contradicts the irreducibility of the representation
(\ref{repr1}).

If $(m_{\infty},\rho_{\infty})=(0,1/2)$, then
$-2m_{\infty}-2\rho_{\infty}+1=0$ and $tb_1^{\bf m}(t)+b_3^{\bf
m}(t)\equiv c={\rm const}$. Hence $u_{\bf m}(t)=-tb_2^{\bf
m}(t)/c$. Note that $c\ne0$, since in the opposite case for all
$t\in D(t^*)$ the function $b_{12}^{\bf m}(z,t)$ has no zeros and
the Riemann--Hilbert problem for scalar Fuchsian equations under
consideration is soluble without apparent singularities, and
$\gamma_{\max}(\chi^*)=2$ (but this contradicts the above
construction).

Thus, if $(m_{\infty},\rho_{\infty})\ne(0,1/2)$, then the poles of
the function $u_{\bf m}(t)$ can be simple only, and if
$(m_{\infty},\rho_{\infty})=(0,1/2)$, then pole orders of the
function $u_{\bf m}(t)$ do not exceed two.
\medskip

{\bf Remark 2.} Alongside formulae for the transition from a
two-dimensional Schlesinger isomono\-dromic family with
$sl(2,{\mathbb C})$-residues to an equation ${\rm P_{VI}}$, there
also exist formulae for the inverse transition (see \cite{JM} or
\cite{Boa}).

Hence, the latter reasonings prove the well known statement about
movable poles of the equation ${\rm P_{VI}}(\alpha,\beta,\gamma,
\delta)$. In the case $\alpha\ne0$ they can be simple only, and in
the case $\alpha=0$ their orders do not exceed two or
$u(t)\equiv\infty$ (see, for instance, \cite{GL}, Ch. VI, \S6).

Indeed, if a solution $u(t)$ of the equation (\ref{PVI})
corresponds to a two-dimensional Schlesinger isomonodromic family
with {\it irreducible} monodromy, then the statement follows from
the above construction ($\alpha\ne0\Longrightarrow
(m_{\infty},\rho_{\infty})\ne(0,1/2)$; $\alpha=0\Longrightarrow
(m_{\infty},\rho_{\infty})=(0,1/2)$, furthermore the case
$\alpha=0$, $u(t)\equiv\infty$ is possible, if the monodromy is
$1$-smaller). If the monodromy of the corresponding family is {\it
reducible}, then $u(t)$ satisfies a Riccati equation (as shown by
M.\,Mazzocco \cite{Maz}), whose movable poles are simple.
\newpage

{\bf\S5. The Riemann--Hilbert problem and Garnier systems}
\medskip

The arguments given above can be extended to general case of $n+3$
singular points $a_1,\ldots,a_n$, $a_{n+1}=0$, $a_{n+2}=1$,
$a_{n+3}=\infty$ and an {\it irreducible non-smaller}
representation
\begin{eqnarray}\label{repr2}
\chi^*: \pi_1({\mathbb C}\setminus
\{a_1,\ldots,a_n,0,1\})\longrightarrow{\rm GL}(2,\mathbb C),
\end{eqnarray}
$a=(a_1,\ldots,a_n)\in D(a^*)$, where $D(a^*)$ is a disk of small
radius centered at the point $a^*$ of the space $({\mathbb
C}\setminus\{0,1\})^n\setminus\bigcup_{i\ne j} \{a_i=a_j\}$.

(Continuing investigations of R.\,Fuchs) R.\,Garnier \cite{Gar}
obtained for $n>1$ the system of non-linear partial differential
equations of second order that must be satisfied by apparent
singularities $\lambda_1(a),\ldots,\lambda_n(a)$ of some Fuchsian
differential equation of second order with singular points
$a_1,\ldots,a_n,0,1,\infty$ and ${\rm SL}(2,\mathbb C)$-monodromy
not depending on the parameter $a$. We supplement these results by
the following reasonings.
\medskip

{\bf Lemma 2.} {\it One has $\gamma_{\max}(\chi^*)=1$ for almost
all $a\in D(a^*)$.}
\smallskip

{\bf Proof.} Consider an arbitrary bundle $F^{\Lambda}$ from the
family $\cal F$ constructed by the representation $\chi^*$. It is
sufficient to prove that $\gamma(F^{\Lambda})\leqslant1$ for all
$a\in D(a^*)$, may be, with the exception of an analytic subset of
codimension one.

If $F^{\Lambda}\cong{\cal O}(k_1)\oplus{\cal O}(k_2)$,
$k_1-k_2>1$, for some $a^0\in D(a^*)$, then we can apply
Bolibrukh's procedure (which was used in the proof of Proposition
1) to get a Schlesinger isomonodromic family of the form
(\ref{addschl}) with an isomonodromic fundamental matrix $Y(z,a)$
of the form (\ref{addschlfund}), where
$$
K={\rm diag}(k'_1,k'_2),\qquad k'_1-k'_2\leqslant1.
$$
This family is defined in the exterior of some analytic subset
$\Theta_{\Lambda}\subset D(a^*)$ of codimension one.

In view of the form of the matrix $Y(z,a)$, the transformation
$y'=z^{-K}y$ transforms this family into the family that is
holomorphic at the infinity. Hence, due to Proposition 2, the
matrix $K$ defines the splitting type of the bundle $F^{\Lambda}$
for $a\in D(a^*)\setminus\Theta_{\Lambda}$ (and
$\gamma(F^{\Lambda})\leqslant1$ for these values of $a$). {\hfill
$\Box$}
\medskip

Thus, in view of the formula (\ref{numbersing}), for almost all
$a\in D(a^*)$ the set of points $a_1,\ldots,a_n, 0, 1, \infty$ and
representation $\chi^*$ can be realized by a Fuchsian differential
equation of second order with $n$ apparent singularities
$u_1(a),\ldots,u_n(a)$. Let us recall how they are connected with
a Garnier system in the case when $\chi^*$ is an ${\rm
SL}(2,\mathbb C)$-representation.

Applying again Corollary 1, let us choose a value of the parameter
$a=a^0=(a_1^0,\ldots,a_n^0)\in D(a^*)$ for which the
representation $\chi^*$ is realized by Fuchsian systems
\begin{eqnarray}\label{fuchs2n}
\frac{dy}{dz}=\left(\sum_{i=1}^{n+2}\frac{B_i^{\bf
m}}{z-a_i^0}\right)y,\qquad {\bf m}=
(m_1,\ldots,m_{n+2},m_{\infty})\in{\mathbb Z}_+^{n+3},
\end{eqnarray}
with the singular points $a_1^0,\ldots,a_n^0$, $a_{n+1}^0=0$,
$a_{n+2}^0=1$, $a_{n+3}^0=\infty$ (here the eigenvalues of the
matrices $B_i^{\bf m}$ are $\pm(m_i+\rho_i)$, and the matrices
$B_{\infty}^{\bf m}=-\sum_{i=1}^{n+2}B_i^{\bf m}$ are diagonal).

Every system of the form (\ref{fuchs2n}) can be included into the
Schlesinger isomonodromic family
\begin{eqnarray}\label{schl2n}
\frac{dy}{dz}=\left(\sum_{i=1}^{n+2}\frac{B_i^{\bf
m}(a)}{z-a_i}\right)y, \qquad B_i^{\bf m}(a^0)=B_i^{\bf m},
\end{eqnarray}
of Fuchsian systems with singularities $a_1,\ldots,a_n,0,1,\infty$
which depends holomorphically on the parameter
$a=(a_1,\ldots,a_n)\in D(a^0)$, furthermore
$\sum_{i=1}^{n+2}B_i^{\bf m}(a)=-B_{\infty}^{\bf m}={\rm
diag}(-m_{\infty}-\rho_{\infty}, m_{\infty}+\rho_{\infty})$.

By Malgrange's theorem the matrix functions
$$
B_i^{\bf m}(a)=\left(\begin{array}{cc} c_i^{\bf m}(a) & b_i^{\bf m}(a) \\
                               d_i^{\bf m}(a) & -c_i^{\bf m}(a)
               \end{array}\right)
$$
can be extended to the universal covering $Z$ of the space
$({\mathbb C}\setminus\{0,1\})^n\setminus\bigcup_{i\ne
j}\{a_i=a_j\}$ as meromorphic functions (holomorphic off the
analytic subset $\Theta_{\bf m}$ of codimension one).

Denote by $B_{\bf m}(z,a)$ the coefficient matrix of the family
(\ref{schl2n}). Since the upper-right element of the matrix
$B_{\infty}^{\bf m}$ equals zero, for every fixed $a$ the same
element of the matrix $z(z-1)(z-a_1)\ldots(z-a_n)B_{\bf m}(z,a)$
is a polynomial $P_{\bf m}(z,a)$ of degree $n$ in $z$. We denote
by $u^{\bf m}_1(a),\ldots,u^{\bf m}_n(a)$ the roots of this
polynomial and define the functions $v^{\bf m}_1(a),\ldots,v^{\bf
m}_n(a)$:
$$
v^{\bf m}_j(a)=\sum_{i=1}^{n+2}\frac{c_i^{\bf m}(a)+m_i+\rho_i}
{u^{\bf m}_j(a)-a_i}, \qquad j=1,\ldots,n.
$$
Then the following statement takes place: {\it the pair $(u^{\bf
m},v^{\bf m}) =(u^{\bf m}_1,\ldots,u^{\bf m}_n, v^{\bf
m}_1,\ldots,v^{\bf m}_n)$ satisfies the Garnier system
$(\ref{garnier})$ with the parameters
$2m_1+2\rho_1,\ldots,2m_{n+2}+2\rho_{n+2},2m_{\infty}+2\rho_{\infty}-1$}
(see \cite{GP}, Cor. 6.2.2).

Thus, using arguments analogous to those given in the case $n=1$,
we get the following statement.
\medskip

{\bf Theorem 4.}{\it The set of the points $a_1,\ldots,a_n$,
$0,1,\infty$ and any irreducible non-smaller ${\rm SL}(2,{\mathbb
C})$-representation $(\ref{repr2})$ can be realized by the family
$($depending on the parameter ${\bf m}\in{\mathbb Z}_+^{n+3})$ of
scalar Fuchsian equations
$$
\frac{d^2w}{dz^2}+p_{\bf m}(z,a)\frac{dw}{dz}+q_{\bf m}(z,a)w=0
$$
with $n$ apparent singularities $($the apparent singular points
$u^{\bf m}_1(a),\ldots,u^{\bf m}_n(a)$ of every equation from this
family and the functions $v^{\bf m}_1(a)={\rm res}\,q_{\bf
m}(z,a)|_{z=u^{\bf m}_1},\ldots,v^{\bf m}_n(a)={\rm res}\,q_{\bf
m}(z,a)|_{z=u^{\bf m}_n}$, $a\in D(a^*)$, form a solution $(u^{\bf
m}(a), v^{\bf m}(a))$ of the Garnier system $(\ref{garnier})$ with
the parameters $2m_1+2\rho_1,\ldots,2m_{n+2}+2\rho_{n+2},
2m_{\infty}+2\rho_{\infty}-1)$}.
\medskip

{\bf Remark 3.} Earlier M.\,Ohtsuki \cite{Oh} has obtained that
the representation (\ref{repr2}) can be realized by a Fuchsian
equation with {\it at most} $n$ apparent singularities (he also
required one of the generating matrices $G_i$ to be
diagonalisable). Here, using Bolibrukh's formula
(\ref{numbersing}), we show that the number of apparent
singularities is $n$ exactly (for almost all locations of
singularities $a_1,\ldots,a_n$).
\smallskip

One can express the coefficients of the polynomial $P_{\bf
m}(z,a)$ in terms of the upper-right elements $b^{\bf m}_i(a)$ of
the matrices $B_i^{\bf m}(a)$. Let
$$
\sigma_1(a)=\sum_{i=1}^{n+2}a_i,\quad \sigma_2(a)=\sum_{1\leqslant
i<j\leqslant n+2}a_ia_j,\quad\ldots,\quad
\sigma_{n+1}(a)=a_1\ldots a_n
$$
be the elementary symmetric polynomials in $a_1,\ldots,a_n$,
$a_{n+1}=0$, $a_{n+2}=1$, and $Q(z)=\prod_{i=1}^{n+2}(z-a_i)$.
Then
$$
P_{\bf m}(z,a)=\sum_{i=1}^{n+2}b^{\bf m}_i(a)\frac{Q(z)}{z-a_i}=:
b_{\bf m}(a)z^n+f_1^{\bf m}(a)z^{n-1}+\ldots+f_n^{\bf m}(a)
$$
(recall that $\sum_{i=1}^{n+2}b^{\bf m}_i(a)=0$). By the Vi\`ete
theorem one has
\begin{eqnarray*}
b_{\bf m}(a)&=&\sum_{i=1}^{n+2}b^{\bf m}_i(a)(-\sigma_1(a)+a_i)=
\sum_{i=1}^{n+2}b^{\bf m}_i(a)a_i=\sum_{i=1}^nb^{\bf m}_i(a)a_i+
b^{\bf m}_{n+2}(a),\\
f_1^{\bf m}(a)&=&\sum_{i=1}^{n+2}b^{\bf m}_i(a)\Bigl(\sigma_2(a)-
\sum_{j=1,j\ne i}^{n+2}a_ia_j\Bigr)=-\sum_{1\leqslant i<j\leqslant
n+2}(b^{\bf m}_i(a)+b^{\bf m}_j(a))a_ia_j.
\end{eqnarray*}

In the similar way,
$$
f_k^{\bf m}(a)=(-1)^k\sum_{1\leqslant i_1<\ldots<i_{k+1}\leqslant
n+2}(b^{\bf m}_{i_1}(a)+\ldots+b^{\bf m}_{i_{k+1}}(a))
a_{i_1}\ldots a_{i_{k+1}}
$$
for each $k=1,\ldots,n$.

It immediately follows from the above formulae and Malgrange's
theorem that {\it the elementary symmetric polynomials
$\sigma_k(u^{\bf m}_1,\ldots,u^{\bf m}_n)=(-1)^kf_k^{\bf
m}(a)/b_{\bf m}(a)$, depending on solutions of the Garnier system
extended to $Z$, are meromorphic functions.}

For $n>1$ a Garnier system generically does not satisfy the
Painlev\'e property (coordinates $(u_1,\ldots,u_n)$ are defined as
roots of a polynomial of degree $n$), but it can be transformed by
a certain (symplectic) transformation $(u,v,a,H)\mapsto(q,p,s,K)$,
$\sum_{i=1}^n(p_idq_i-K_ids_i)=\sum_{i=1}^n(v_idu_i-H_ida_i)$,
into a Hamiltonian system satisfying the Painlev\'e property (see
\cite{GP}, Ch. III, \S7).

By Theorem 1 for each function $f_k^{\bf m}(a)$ extended to $Z$
and any point $a^*$ of the $\Theta$-divisor of the family
(\ref{schl2n}) one has $\Sigma_{a^*}(f_k^{\bf m})\geqslant-n-1$.
Similarly to the case $n=1$, here we can tell something about the
behaviour of the function $b_{\bf m}(a)$ along $\Theta_{\bf m}$.
\medskip

{\bf Lemma 3.} {\it Consider the family (\ref{schl2n}) with the
irreducible non-smaller monodromy $\chi^*$, and the function
$b_{\bf m}(a)$ constructed by the residue matrices $B_i^{\bf
m}(a)$.

{\rm i)} In the case $(m_{\infty},\rho_{\infty})=(0,1/2)$ one has
$b_{\bf m}(a)\equiv{\rm const}\ne0$;

{\rm ii)} In the case $(m_{\infty},\rho_{\infty})\ne(0,1/2)$ the
set $\{a\in Z\,|\,b_{\bf m}(a)=0\}$ is an analytic submanifold of
codimension one in $Z$, and if the function $b_{\bf m}(a)$ is
holomorphic at a point $a^0\in Z$, so are the functions $f_k^{\bf
m}(a)$.}
\smallskip

{\bf Proof.} By Lemma 1 we have $db_{\bf m}(a)=
(-2m_{\infty}-2\rho_{\infty}+1)\sum_{i=1}^nb_i^{\bf m}(a)da_i$.

i) In the case $(m_{\infty},\rho_{\infty})=(0,1/2)$ one has
$db_{\bf m}(a)\equiv0$ for all $a\in D(a^*)$, hence $b_{\bf
m}(a)\equiv{\rm const}\ne0$. Indeed, if $b_{\bf m}(a)\equiv0$,
then $P_{\bf m}(z,a)$ is a polynomial of degree $n-1$ in $z$.
Therefore, for every $a\in D(a^*)$ the representation $\chi^*$ is
realized by a scalar Fuchsian equation with at most $n-1$ apparent
singularities (which are the roots of $P_{\bf m}(z,a)$) and
$\gamma_{\max}(\chi^*)>1$, which contradicts Lemma 2.

ii) In the case $(m_{\infty},\rho_{\infty})\ne(0,1/2)$ one has
$\theta=-2m_{\infty}-2\rho_{\infty}+1\ne 0$, and
\begin{eqnarray}\label{formulae}
b_i^{\bf m}(a)&=&\frac1{\theta}\frac{\partial b_{\bf
m}(a)}{\partial a_i},\qquad i=1,\ldots,n;\nonumber\\
b_{n+2}^{\bf m}(a)&=&b_{\bf m}(a)-\sum_{i=1}^nb_i^{\bf m}(a)a_i,
\quad b_{n+1}^{\bf m}(a)=-b_{n+2}^{\bf m}(a)-\sum_{i=1}^nb_i^{\bf
m}(a).
\end{eqnarray}
Thus, if the function $b_{\bf m}(a)$ is holomorphic at some point
$a^0\in Z$, so are the functions $b_i^{\bf m}(a)$,
$i=1,\ldots,n+2$, and hence, the functions $f_k^{\bf m}(a)$.

If for some $a^0\in\{b_{\bf m}(a)=0\}$ one has $db_{\bf
m}(a^0)\equiv0$, then $\sum_{i=1}^nb_i^{\bf m}(a^0)da_i\equiv0$
and $b^{\bf m}_1(a^0)=\ldots=b^{\bf m}_n(a^0)=0$. Taking into
consideration the relations (\ref{formulae}), one gets also
$b^{\bf m}_{n+2}(a^0)=0$ and $b^{\bf m}_{n+1}(a^0)=0$. This
contradicts the irreducibility of the representation $\chi^*$.
{\hfill $\Box$}
\medskip

As a consequence of Theorem 1 and Lemma 3, one gets the following
statement.
\medskip

{\bf Theorem 5.} {\it Denote by $\Delta_i$ the polar loci of the
functions $\sigma_i(u^{\bf m}_1(a),\ldots,u^{\bf m}_n(a))$
extended to $Z$, respectively $($in the conditions of Theorem
4$)$. Then

{\rm a)} in the case $(m_{\infty},\rho_{\infty})=(0,1/2)$ one has
$\Sigma_{a^*}(\sigma_i)\geqslant-n-1$ for any point
$a^*\in\Delta_i$;

{\rm b)} in the case $(m_{\infty},\rho_{\infty})\ne(0,1/2)$ one
has $\Sigma_{a^*}(\sigma_i)\geqslant-n$ for any point
$a^*\in\Delta_i\setminus\Delta^0$, where $\Delta^0\subset\Delta_i$
is some subset of positive codimension $($or the empty set$)$;

{\rm c)} in the case $(m_{\infty},\rho_{\infty})\ne(0,1/2)$,
$a^*\in\Delta^0$, one can estimate the order $\kappa$ of each
irreducible component of $\Delta_i\cap D(a^*)$ as follows:
$\kappa\geqslant-n$.}

{\bf Proof.} Recall that $\sigma_i(u^{\bf m}_1,\ldots,u^{\bf
m}_n)=(-1)^if_i^{\bf m}(a)/b_{\bf m}(a)$ and
$\Sigma_{a^*}(f_i^{\bf m})\geqslant-n-1$ for any point $a^*$ of
the $\Theta$-divisor of the family (\ref{schl2n}).

Therefore, the statement a) of the theorem is a consequence of
Lemma 3, i).

b) As follows from Lemma 3, ii), the points $a^*\in\Delta_i$ can
be of two types: such that $b_{\bf m}(a^*)=0$ (then
$\Sigma_{a^*}(\sigma_i)\geqslant-1$) or that belong to
$\Theta_{\bf m}$.

Denote by $\Delta^0$ the set of the points of $\Theta_{\bf m}$
that are ambiguous for $b_{\bf m}(a)$. Then in a neighbourhood of
any point $a^*\in\Theta_{\bf m}\setminus\Delta^0$ each function
$f_i^{\bf m}(a)$ can be presented in the form
\begin{eqnarray}\label{divisor1}
f_i^{\bf m}(a)=\frac{g(a)}{\tau_1^{k_1}(a)\ldots\tau_r^{k_r}(a)},
\quad k_1+\ldots+k_r\leqslant n+1,
\end{eqnarray}
where $\tau_i(a)$, $g(a)$ are holomorphic near $a^*$, furthermore
$\tau_i(a)$ are irreducible at $a^*$, just as
\begin{eqnarray}\label{divisor2}
b_{\bf m}(a)=\frac{h(a)}{\tau_1^{j_1}(a)\ldots\tau_r^{j_r}(a)},
\quad j_1+\ldots+j_r\geqslant1,
\end{eqnarray}
where $h(a)$ is holomorphic near $a^*$, $h(a^*)\ne0$. Thus,
$$
\frac{f_i^{\bf m}(a)}{b_{\bf m}(a)}=
\frac{g(a)}{\tau_1^{k_1}(a)\ldots\tau_r^{k_r}(a)}:
\frac{h(a)}{\tau_1^{j_1}(a)\ldots\tau_r^{j_r}(a)}=
\frac{g(a)/h(a)}{\tau_1^{k_1-j_1}(a)\ldots\tau_r^{k_r-j_r}(a)},
$$
therefore,
$$
\Sigma_{a^*}(\sigma_i)=-(k_1-j_1)-\ldots-(k_r-j_r)\geqslant-n.
$$

c) In a neighbourhood of a point $a^*\in\Delta^0$ the
decompositions (\ref{divisor1}), (\ref{divisor2}) take place for
the functions $f_i^{\bf m}(a)$, $b_{\bf m}(a)$ respectively, but
$h(a^*)=0$. However, due to Lemma 3, ii), all irreducible factors
of $h(a)$ in its decomposition $h(a)=h_1(a)\ldots h_s(a)$ near
$a^*$ are distinct (we can assume also that none of $h_i$
coincides with some of $\tau_j$). One also has $k_i=0$, if $j_i=0$
($b_{\bf m}(a)$ is holomorphic along
$\{\tau_i(a)=0\}\Longrightarrow f_i^{\bf m}(a)$ is holomorphic
along $\{\tau_i(a)=0\}$). Therefore, $k_i-j_i\leqslant n$, and the
statement c) follows from the decomposition
$$
\frac{f_i^{\bf m}(a)}{b_{\bf m}(a)}=\frac{g(a)}{h_1(a)\ldots
h_s(a)\,\tau_1^{k_1-j_1}(a)\ldots\tau_r^{k_r-j_r}(a)}.
$$
{\hfill $\Box$}
\medskip

Alongside formulae for the transition from a two-dimensional
Schlesinger isomonodromic family with $sl(2,{\mathbb C})$-residues
to a Garnier system, there also exist formulae for the inverse
transition (see \cite{GP}, Ch. III, \S6.3). Hence, the latter
theorem implies some addition to Garnier's theorem \cite{Gar}
(which claims that the elementary symmetric polynomials of
solutions of a Garnier system are meromorphic on $Z$).
\medskip

{\bf Theorem 5 bis.} {\it Consider a solution $(u(a),v(a))$ of the
Garnier system $(\ref{garnier})$, that corresponds to a
two-dimensional Schlesinger isomonodromic family with irreducible
monodromy, and the polar loci $\Delta_i$ of the functions
$\sigma_i(u_1(a),\ldots,u_n(a))$ meromorphic on $Z$. Then

{\rm a)} in the case $\theta_{\infty}=0$ and the non-smaller
monodromy one has $\Sigma_{a^*}(\sigma_i)\geqslant-n-1$ for any
point $a^*\in\Delta_i$;

{\rm b)} in the case $\theta_{\infty}\ne0$ one has
$\Sigma_{a^*}(\sigma_i)\geqslant-n$ for any point
$a^*\in\Delta_i$, may be, with the exception of some subset
$\Delta^0\subset\Delta_i$ of positive codimension $($for whose
points $a^0$ we have the estimate $\kappa\geqslant-n$ for the
order $\kappa$ of each irreducible component of $\Delta_i\cap
D(a^0)\,)$.}
\medskip

{\bf Remark 4.} M.\,Mazzocco \cite{Maz2} has shown that the
solutions of the Garnier system (\ref{garnier}), that correspond
to two-dimensional Schlesinger isomonodromic families with {\it
reducible} monodromy, are classical functions (in each variable,
in sense of Umemura \cite{Um}) and can be expressed via Lauricella
hypergeometric equations (see \cite{GP}, Ch. III, \S9). Thus,
Theorem 5 bis can be applied, for example, to non-classical
solutions of Garnier systems.


\begin{thebibliography}{99}

\bibitem{AB}
D.\,V.\,Anosov, A.\,A.\,Bolibruch, The Riemann--Hilbert problem.
{\it Aspects Math.}, {\bf E 22}, Braunschweig: Vieweg, 1994.

\bibitem{AL}
D.\,V.\,Anosov, V.\,P.\,Leksin, Andrei Andreevich Bolibrukh in
life and science. {\it Russian Math. Surveys}, 2004, {\bf59}(6),
1009--1028.

\bibitem{Boa}
P.\,P.\,Boalch, Some explicit solutions to the Riemann--Hilbert
problem. {\it IRMA Lect. Math. Theor. Phys.}, 2006, {\bf9},
85--112. (available at http://www.arxiv.org)

\bibitem{Bo1}
A.\,A.\,Bolibrukh, The Riemann--Hilbert problem on the complex
projective line (in Russian). {\it Mat. Zametki}, 1989,
{\bf46}(3), 118--120.

\bibitem{Bo4}
A.\,A.\,Bolibruch, Vector bundles associated with monodromies and
asymptotics of Fuchsian systems. {\it J. Dynam. Control Systems},
1995, {\bf1}(2), 229--252.

\bibitem{Bo2}
A.\,A.\,Bolibruch, On isomonodromic deformations of Fuchsian
systems. {\it J. Dynam. Control Systems}, 1997, {\bf3}(4),
589--604.

\bibitem{Bo5}
A.\,A.\,Bolibruch, On orders of movable poles of the Schlesinger
equation. {\it J. Dynam. Control Systems}, 2000, {\bf6}(1),
57--74.

\bibitem{Bo3}
A.\,A.\,Bolibrukh, On the tau-function for the Schlesinger
equation of isomonodromic deformations. {\it Math. Notes}, 2003,
{\bf74}(2), 177--184.

\bibitem{Bo6}
A.\,A.\,Bolibruch, Inverse monodromy problems of the analytic
theory of differential equations. {\it Mathematical events of the
twentieth century}, Berlin: Springer, 2006, 49--74.

\bibitem{Ch}
E.\,M.\,Chirka, Complex analytic sets. {\it Mathematics and its
applications}, {\bf46}, Dordrecht: Kluwer, 1989.

\bibitem{Fu}
L.\,Fuchs, Zur Theorie der linearen Differentialgleichungen mit
ver\"anderlichen Koefficienten. {\it J.~Reine Angew. Math.}, 1868,
{\bf68}, 354--385.

\bibitem{Fu1}
R.\,Fuchs, Sur quelques \'equations diff\'erentielles lin\'eaires
du second ordre. {\it C. R. Acad. Sci. Paris}, 1905, {\bf141},
555--558.

\bibitem{Ga}
B.\,Gambier, Sur les \'equations diff\'erentielles du second ordre
et du premier degr\'e dont l'int\'egrale g\'en\'erale est \`a
points critiques fix\'es. {\it C. R. Acad. Sci. Paris}, 1906,
{\bf142}, 266--269.

\bibitem{Gar}
R.\,Garnier, Sur des \'equations diff\'erentielles du troisi\`eme
ordre dont l'int\'egrale g\'en\'erale est uniforme et sur une
classe d'\'equations nouvelles d'ordre sup\'erieur dont
l'int\'egrale g\'en\'erale a ses points critiques fixes. {\it Ann.
Sci. \'Ecole Norm. Sup.}, 1912, {\bf29}, 1--126.

\bibitem{Go}
R.\,R.\,Gontsov, On solutions of the Schlesinger equation in the
neighbourhood of the Malgrange $\Theta$-divisor. {\it Math.
Notes}, 2008, {\bf83}(5), 707--711.

\bibitem{Go2}
R.\,R.\,Gontsov, Apparent singularities of Fuchsian equations, the
Painlev\'e VI equation, and Garnier systems. {\it Doklady Math.},
2009, {\bf79}(2), 176--179.

\bibitem{GL}
V.\,I.\,Gromak, N.\,A.\,Lukashevich, The analytic solutions of the
Painlev\'e equations (in Russian). Minsk: Universitetskoye, 1990.

\bibitem{Ha}
P.\,Hartman, Ordinary differential equations. New York: Wiley,
1964.

\bibitem{GP}
K.\,Iwasaki, H.\,Kimura, S.\,Shimomura, M.\,Yoshida, From Gauss to
Painlev\'e. {\it Aspects Math.}, {\bf E 16}, Braunschweig: Vieweg,
1991.

\bibitem{JM}
M.\,Jimbo, T.\,Miwa, Monodromy preserving deformations of linear
differential equations with rational coefficients. II. {\it
Physica D}, 1981, {\bf2}(3), 407--448.

\bibitem{Le}
A.\,Levelt, Hypergeometric functions. {\it Proc. Konikl. Nederl.
Acad. Wetensch. Ser.~A}, 1961, {\bf64}, 361--401.

\bibitem{Ma}
B.\,Malgrange, Sur les d\'eformations isomonodromiques. I.
Singularit\'es r\'eguli\`eres. {\it Progr. Math.}, 1983, {\bf37},
401--426.

\bibitem{Maz}
M.\,Mazzocco, Rational solutions of the Painlev\'e VI equation.
{\it J. Phys. A}, 2001, {\bf34}(11), 2281--2294. (available at
http://www.arxiv.org)

\bibitem{Maz2}
M.\,Mazzocco, The geometry of the classical solutions of the
Garnier systems. {\it I.R.M.N.}, 2002, {\bf12}, 613--646.
(available at http://www.arxiv.org)

\bibitem{Oh}
M.\,Ohtsuki, On the number of apparent singularities of a linear
differential equation. {\it Tokyo J. Math.}, 1982, {\bf5}(1),
23--29.

\bibitem{Ok}
K.\,Okamoto, Isomonodromic deformation and Painlev\'e equations
and the Garnier system. {\it J. Fac. Sci. Univ. Tokyo Sec. IA,
Math.}, 1986, {\bf33}, 575--618.

\bibitem{Po}
A.\,Poincar\'e, Sur les groupes des \'equations lin\'eaires. {\it
Acta Math.}, 1884, {\bf4}, 201--311.

\bibitem{Ri}
B.\,Riemann, Zwei allgemeine Lehrs\"atze \"uber lineare
Differentialgleichungen mit algebraichen Koefficienten. {\it Math.
Werke}, 1892, 357--369.

\bibitem{Sch}
L.\,Schlesinger, Uber Losungen gewisser Differentialgleichungen
als Funktionen der singularen Punkte. {\it J.~Reine Angew. Math.},
1905, {\bf129}, 287--294.

\bibitem{Um}
H.\,Umemura, Birational automorphism groups and differential
equations. {\it Nagoya Math. J.}, 1990, {\bf119}, 1--80.

\bibitem{VG}
I.\,V.\,Vyugin, R.\,R.\,Gontsov, Additional parameters in inverse
monodromy problems. {\it Russian Acad. Sci. Sb. Math.}, 2006,
{\bf197}(12), 1753--1773.

\end{thebibliography}
\end{document}